\documentclass[11pt, a4paper]{article}
\usepackage{amsmath}
\usepackage{amsthm}
\usepackage{mathtools}
\usepackage{mathrsfs}
\usepackage[numbers]{natbib}
\usepackage{lmodern}
\usepackage{mathpazo}
\usepackage{amssymb}
\usepackage[T1]{fontenc}
\usepackage{microtype}
\numberwithin{equation}{section}
\usepackage[margin=1in]{geometry}
\usepackage{indentfirst}
\usepackage{enumitem}
\usepackage{array}
\usepackage{placeins}
\usepackage{xcolor}
\usepackage{hyperref}
\hypersetup{
    colorlinks=true,
    linkcolor=black,
    citecolor=black,
    urlcolor=black,
    pdftitle={Reflections of Nichols algebras in braided monoidal categories},
    pdfauthor={Bowen Li},
    pdfsubject={Nichols algebras in braided monoidal categories},
    pdfkeywords={Nichols algebra, braided monoidal category, Yetter--Drinfeld module, reflection}
}
\usepackage{cleveref}
\usepackage{needspace}

\newcommand{\C}{\mathcal C}
\newcommand{\D}{\mathcal D}
\newcommand{\B}{\mathcal B}
\newcommand{\YD}{\mathcal{YD}}
\newcommand{\one}{\mathbf 1}
\DeclareMathOperator{\id}{id}
\DeclareMathOperator{\ad}{ad}
\DeclareMathOperator{\co}{co}
\DeclareMathOperator{\im}{im}
\DeclareMathOperator{\Ind}{Ind}
\newcommand{\rat}{\mathrm{rat}}
\definecolor{CrossRefColor}{RGB}{0,70,160}
\definecolor{CiteColor}{RGB}{0,110,60}
\definecolor{URLColor}{RGB}{130,40,140}

\hypersetup{
  colorlinks=true,
  linkcolor=CrossRefColor,
  citecolor=CiteColor,
  urlcolor=URLColor
}
\newtheorem{theorem}{Theorem}[section]
\newtheorem{proposition}[theorem]{Proposition}
\newtheorem{lemma}[theorem]{Lemma}
\newtheorem{corollary}[theorem]{Corollary}
\newtheorem{definition}[theorem]{Definition}
\newtheorem{remark}[theorem]{Remark}
\setlist[enumerate,1]{label=\textnormal{(\arabic*)},leftmargin=*,itemsep=3pt}

\makeatletter
\newcommand\blfootnote[1]{%
  \begingroup
  \renewcommand\thefootnote{}%
  \footnotetext[0]{%
    \begin{minipage}[t]{\dimexpr\linewidth-1.8em\relax}
      \setlength{\parindent}{0pt}%
      #1
    \end{minipage}%
  }%
  \endgroup
}
\makeatother

\title{Reflections of Nichols algebras in braided monoidal categories}
\author{Bowen Li}
\date{}

\definecolor{ApplicationColor}{HTML}{000000}

\definecolor{red}{HTML}{000000}
\ifdefined\CleanCopy
  \definecolor{RevisionColor}{HTML}{000000}
  \definecolor{ConcisionColor}{HTML}{000000}
  \definecolor{ApplicationColor}{HTML}{000000}
  \definecolor{red}{HTML}{000000}
  \definecolor{blue}{HTML}{000000}
  \hypersetup{hidelinks}
\fi

\begin{document}
\maketitle
 
\blfootnote{
 {\textit{Affiliation:} School} of Mathematical Sciences, Shanghai Jiao Tong University.

\textit{Email:}
\href{mailto:lbwwysusu@gmail.com}{lbwwysusu@gmail.com}
}

\begin{abstract}
Let $\C$ be an essentially small abelian rigid braided monoidal
category. We prove a reflection theorem for Nichols algebras in
$\Ind(\C)$ and show that tuples admitting all reflections define
connected semi-Cartan graphs.

As applications, we prove that a finite Nichols algebra has a generating
object of length at most one in $PSU(2)_k$ for odd $k\ge3$, in
$SU(2)_{4m}$ for $m\ge1$, and in the representation categories of the
Virasoro minimal models $\mathcal M(2,p)$ for odd $p\ge5$.
In all three families, we prove that  every finite
pre-Nichols algebra is a Nichols algebra.
\end{abstract}

\section{Introduction}

{Reflections of Nichols algebras of finite-dimensional
semisimple Yetter--Drinfeld modules over Hopf algebras were
constructed by Andruskiewitsch, Heckenberger and
Schneider~\cite{AHS10}. Heckenberger and Schneider developed
the associated root systems and Weyl groupoids~\cite{HS10},
and a construction of reflections using dually paired Hopf
algebras~\cite{HSdualpair}, see also~\cite{rootsys}.}
Wolf~\cite{Wolf21} studied reflection functors
for Yetter--Drinfeld modules over Nichols systems in this setting,
see also~\cite{Wolf24}.
 Reflection theory for coquasi-Hopf algebras was established in \cite{reflection2} and \cite{reflection3}.

Barvels, Lentner and Schweigert~\cite{Simon} constructed partial
dualization in general braided monoidal categories. Their construction
uses a Hopf pairing with an inverse copairing for the whole Hopf
algebra {~\cite[Definition~2.15]{Simon}}.
For graded Hopf algebras with infinitely many nonzero components,
degreewise duality and rationality provide the appropriate replacement.
This approach is already present in the Hopf and coquasi-Hopf
settings {~\cite{HSdualpair,reflection2}}. The question considered here is whether the coinvariant
identification and simplicity results needed for reflections hold
without such a realization of the ambient category.

Lentner~\cite[Problem~4, p.~28]{Lentner26} asks for a proof that reflected
components remain irreducible in a general categorical setting with
duals. He also suggests that the Nichols algebra identification of
coinvariants extends to general braided tensor
categories~\cite[Remark~5.2.2, p.~64]{Lentner26}.
We prove these assertions in an essentially small abelian rigid braided
monoidal category $\C$. Nichols algebras are formed in
$\D=\Ind(\C)$.

Let $M=(M_1,\ldots,M_{ {r}})$ be a tuple of nonzero simple objects of
$\C$, and write $\B(M)=\B(\bigoplus_{ {j}=1}^{{ {r}}}M_j)$.
 {Fix $i$. Radford's projection theorem
\cite[Theorem~3.10]{Simon} equips $\B(M)^{\co\B(M_i)}$ with a
Hopf algebra structure in ${}_{\B(M_i)}^{\B(M_i)}\YD(\D)$,
whose action is the restriction of the adjoint action.
Theorem~\ref{coinv:nichols} identifies this Hopf algebra as a
Nichols algebra:}
\[
 \B(M)^{\co\B(M_i)} \cong
 \mathcal{B}((\operatorname{ad}\B(M_i))(\bigoplus_{j\neq i}M_j)).
\]
{Moreover, each $\ad(\B(M_i))(M_j)$, for $j\ne i$,
is simple in ${}_{\B(M_i)}^{\B(M_i)}\YD(\D)$
by Proposition~\ref{coinv:simple}.}
Now suppose that
$(\ad M_i)^{ {n}}(M_j)=0$ for all sufficiently large ${ {n}}$, for each $j\ne i$.
Let $m_{ij}$ be the largest ${ {n}}$ for which
$(\ad M_i)^{ {n}}(M_j)\ne0$, and set
\[
 R_i(M)_i=M_i^*,\qquad
 R_i(M)_j=(\ad M_i)^{m_{ij}}(M_j)\quad(j\ne i).
\]
 {Then we say that $M$ admits the $i$-th reflection.}
We prove that $R_i(M)_j$ are simple  for all $1\leq j \leq { {r}}$ by
Lemmas~\ref{pre:simple-dual} and~\ref{refl:rational}.
Theorem~\ref{refl:main} gives an isomorphism of Hopf algebras in $\D$
\[
 \B(R_i(M))\cong
 \Omega_i\bigl(\B(M)^{\co\B(M_i)}\bigr)\#\B(M_i^*).
\]
Here $\Omega_i$ is the braided equivalence between
the rational Yetter--Drinfeld categories over $\B(M_i)$ and
$\B(M_i^*)$ given by Theorem~\ref{dual:omega}.
 We say that $M$ admits all reflections if every finite
sequence of reflections is defined. Such a tuple determines a
connected semi-Cartan graph, see Remark~\ref{refl:semi-cartan}.
  
The reflection theorem gives a uniform restriction on the
generators of finite Nichols algebras in several families of
braided fusion categories. For $PSU(2)_k$ with odd $k\ge3$,
$SU(2)_{4m}$ with $m\ge1$, and
$\operatorname{Rep}\mathcal M(2,p)$ with odd $p\ge5$, we prove
\[
 V\in\C,\qquad \B(V)\text{ finite}
 \quad\Longrightarrow\quad \operatorname{length}(V)\le1
\]
(Theorems~\ref{app:psu} and~\ref{app:su2}, and
Corollary~\ref{app:minimal}).
Here finite means that the underlying object belongs to $\C$.
Thus every nontrivial finite Nichols algebra in these categories
must have a simple generating object.

These families arise from quantum groups and Virasoro minimal
models \cite{RSW09,GMR24,BD12,Huang}, and illustrate the scope
of an intrinsic reflection theorem.
Each contains a simple object of noninteger Frobenius--Perron
dimension. Hence none is tensor equivalent to the representation
category of a finite-dimensional quasi-Hopf algebra
\cite[Proposition~6.1.14]{tensor}.
The restriction also holds in nonunitary categories, including
the Yang--Lee category $\operatorname{Rep}\mathcal M(2,5)$.

The restriction on generators yields a generation theorem for
finite Hopf algebras in these categories.
Combined with Meir's results
\cite[Theorems~1.4 and~1.6]{defofnichols}, it implies that every
finite pre-Nichols algebra is a Nichols algebra.
Moreover, every Hopf algebra $H\in\C$ satisfying
$H/J(H)\simeq\one$ and $\operatorname{corad}(H)\simeq\one$
is isomorphic to $\B(P(H))$, with $P(H)$ zero or simple
(Corollary~\ref{app:psu-hopf}).

{Section~\ref{sec:preliminaries} collects the categorical preliminaries.}
Sections~\ref{sec:duality}--\ref{sec:reflection} develop rational
Yetter--Drinfeld duality, identify the coinvariants as an internal
Nichols algebra, and establish the simplicity and reflection results.
Section~\ref{sec:examples} applies these results to finite Nichols
and Hopf algebras in braided fusion categories.

\section{Preliminaries}\label{sec:preliminaries}

\subsection{The ambient category}\label{pre:ambient}

Let $\C$ be an essentially small abelian {rigid braided monoidal category},
with unit object $\one$ and braiding $c$.
{By braided monoidal coherence theorem, we may  assume that
$\C$ is strict.} We work in
$\D=\Ind(\C)$ and identify $\C$ with its canonical image in $\D$.

The category $\D$ is a Grothendieck abelian category, and the embedding
$\C\longrightarrow\D$ is fully faithful and exact
\cite[Theorem~8.6.5]{KS06}.
We identify subobjects whose inclusions are related by an
isomorphism over the ambient object. For $Z\hookrightarrow X$ and
$f:X\to Y$, we write $f(Z)$ for the image of the composite
$Z\to X\to Y$.

\begin{lemma}\label{pre:finite-presentability}
\begin{enumerate}
 \item Every object $X\in\C$ is finitely presentable in $\D$,
 that is, $\operatorname{Hom}_{\D}(X,-)$ preserves filtered colimits.
 In particular, every morphism $X\to\varinjlim_iY_i$  factors through some $Y_i$.
 \item For every family $(Y_\lambda)_{\lambda\in\Lambda}$ in $\D$,
let
$\pi_\lambda:\bigoplus_{\mu\in\Lambda}Y_\mu\longrightarrow Y_\lambda
$
be the canonical projection.
For any $X\in\D$ and morphisms
\[
 f,g:X\longrightarrow\bigoplus_{\lambda\in\Lambda}Y_\lambda,
\]
one has $f=g$ if and only if
$\pi_\lambda f=\pi_\lambda g$ for every $\lambda\in\Lambda$.
\end{enumerate}
\end{lemma}

\begin{proof}

\textup{(1)} Realize $\D$ as the full subcategory of
$\operatorname{Fun}(\C^{\mathrm{op}},\mathbf{Set})$ consisting of
filtered colimits of representable functors. Under this realization,
$X\in\C$ is identified with $\operatorname{Hom}_{\C}(-,X)$.
Thus the Yoneda lemma gives
\[
 \operatorname{Hom}_{\D}(X,Y)\cong Y(X)
 \qquad(Y\in\D).
\]
By \cite[Theorem~6.1.8]{KS06}, filtered colimits in $\D$ are
comleted pointwise in this functor category. Hence
\[
\begin{aligned}
 \operatorname{Hom}_{\D}\bigl(X,\varinjlim_i Y_i\bigr)
 \cong \bigl(\varinjlim_i Y_i\bigr)(X)
\cong \varinjlim_i Y_i(X)
 \cong \varinjlim_i\operatorname{Hom}_{\D}(X,Y_i).
\end{aligned}
\]
Every element of the last colimit is represented by a morphism
$X\to Y_i$ for some $i$, which gives the asserted factorization.

(2)
Only the reverse implication needs proof.
Set $h=f-g$. Then $\pi_\lambda h=0$ for every $\lambda$.
First suppose that $X\in\C$.
The direct sum is the filtered colimit of its finite subsums:
\[
 \bigoplus_{\lambda\in\Lambda}Y_\lambda
 \cong
 \varinjlim_{\substack{F\subseteq\Lambda\\F\text{ finite}}}
 \bigoplus_{\lambda\in F}Y_\lambda.
\]
Since $X$ is finitely presentable by \textup{(1)}, there is
a finite subset $F\subseteq\Lambda$ such that $h$ factors as
\[
 X\xrightarrow{h_F}\bigoplus_{\lambda\in F}Y_\lambda
 \xrightarrow{\iota_F}\bigoplus_{\lambda\in\Lambda}Y_\lambda,
\]
where $\iota_F$ is the canonical inclusion.
For each $\lambda\in F$, the $\lambda$-component of $h_F$
is $\pi_\lambda h=0$.
A finite direct sum is also a product, so $h_F=0$.
Thus $h=0$.

For arbitrary $X\in\D$, write
$X=\varinjlim_i X_i$ with $X_i\in\C$, and let
$u_i:X_i\to X$ be the structure morphisms.
For every $i$ and $\lambda$,
\[
 \pi_\lambda h u_i=0.
\]
The preceding argument gives $h u_i=0$ for every $i$.
The universal property of the colimit then gives $h=0$,
and hence $f=g$.
\end{proof}

The tensor product and the braiding extend to $\D$. The extended
tensor product is exact and preserves colimits in each variable.
Indeed, for $X\in\C$, its left dual $X^*$ remains a left dual in
$\D$. Duality and the braiding give natural isomorphisms
\[
\begin{aligned}
 \operatorname{Hom}_{\D}(X^*\otimes Y,Z)
 &\cong \operatorname{Hom}_{\D}(Y,X\otimes Z),\\
 \operatorname{Hom}_{\D}(X\otimes Y,Z)
 &\cong \operatorname{Hom}_{\D}(Y\otimes X,Z)
 \cong \operatorname{Hom}_{\D}(Y,Z\otimes X^*).
\end{aligned}
\]
Hence $X\otimes-$ has left adjoint $X^*\otimes-$ and right adjoint
$-\otimes X^*$, so it is exact and preserves colimits.
For an arbitrary object $X=\varinjlim_i X_i$ of $\D$, expressed
as a filtered colimit with $X_i\in\C$, we have
\[
 X\otimes-\cong\varinjlim_i(X_i\otimes-).
\]
This functor preserves colimits because colimits commute with
colimits, and is exact because filtered colimits in $\D$ are exact.
The braiding gives the same conclusions for $-\otimes X$.

\begin{lemma}\label{pre:simple-dual}
If $X$ is a simple object of $\C$, then $X^*$ is simple
in $\C$.
\end{lemma}

\begin{proof}
Left and right duality are mutually inverse contravariant
equivalences. Hence they are exact and preserve nonzero objects.
In particular, $X^*\neq0$.

Let $Y\subseteq X^*$ be a nonzero subobject and let $Z=X^*/Y$.
Applying the right dual functor ${^*(-)}$ gives an exact sequence
\[
 0\longrightarrow {^*Z}\longrightarrow X
 \longrightarrow {^*Y}\longrightarrow0,
\]
where we use the canonical isomorphism ${^*(X^*)}\simeq X$.
Since $X$ is simple and ${^*Y}\neq0$, we have ${^*Z}=0$.
Thus $Z=0$, so the inclusion $Y\hookrightarrow X^*$ is an
isomorphism.
\end{proof}

\begin{lemma}\label{pre:simple-ind}
Every simple object of $\C$ is simple in $\D$.
\end{lemma}

\begin{proof}
Let $L$ be simple in $\C$, and let $X\longrightarrow L$ be a nonzero
subobject in $\D$. Write $X=\varinjlim_i X_i$ with $X_i\in\C$.
For some $i$, the composite $X_i\longrightarrow X\longrightarrow L$
is nonzero. Its image, comleted in $\C$, is $L$.
Thus $X\longrightarrow L$ is an epimorphism, and hence an isomorphism.
\end{proof}

An $\mathbb N_0$-grading on an object $X\in\D$ is a family
$(X(n))_{n\geq0}$ together with a specified isomorphism
\[
 \bigoplus_{n\geq0}X(n)\cong X.
\]
These direct sums exist in $\D$. We identify each $X(n)$
with the corresponding direct summand of $X$.
The grading is called locally finite if $X(n)\in\C$
for every $n\geq0$. 
Direct sums in $\D$ are exact. Thus, if a morphism preserves
degree, its kernel, image, and cokernel are the direct sums
of those of its homogeneous components.

Since the tensor product preserves direct sums, the tensor
product of graded objects has the grading
\[
 (X\otimes Y)(n)
 =\bigoplus_{p+q=n}X(p)\otimes Y(q).
\]
A graded bialgebra is a bialgebra
$H=\bigoplus_{n\geq0}H(n)$ whose structure morphisms preserve
degree, with $\one$ concentrated in degree zero.
Thus multiplication and comultiplication restrict to morphisms
$
 H(p)\otimes H(q)\longrightarrow H(p+q),
\
 H(n)\longrightarrow\bigoplus_{p+q=n}H(p)\otimes H(q),
$
the unit factors through $H(0)$, and the counit vanishes
on $H(n)$ for $n>0$.
A graded Hopf algebra is a graded bialgebra with a
degree-preserving antipode.
We call a graded bialgebra connected if its unit induces
an isomorphism $\one\cong H(0)$.

Every connected graded bialgebra $H$ has an antipode, defined by
\[
 S_H|_{H(n)}
 =\sum_{r=0}^{n}(\eta_H\varepsilon_H-\id_H)^{*r}|_{H(n)}.
\]
Here $*$ denotes convolution, and the zeroth power is
$\eta_H\varepsilon_H$.
The morphism $\eta_H\varepsilon_H-\id_H$ vanishes on $H(0)$.
 {A direct comletation gives}
$
 S_H*\id_H=\eta_H\varepsilon_H=\id_H*S_H.
$
The primitive subobject of a bialgebra $H$ in $\D$ is
\[
 P(H)=\ker\bigl(\Delta_H-\id_H\otimes\eta_H
                         -\eta_H\otimes\id_H\bigr).
\]
If $H$ is graded, then $P(H)$ is a graded subobject. If $H$ is
connected, then $P(H)(0)=0$ and $H(1)\subseteq P(H)$.

\subsection{Nichols algebras}\label{pre:nichols}

For $V\in\D$, let
\[
 T(V)=\bigoplus_{n\geq0}V^{\otimes n},\qquad V^{\otimes0}=\one.
\]
The tensor algebra $T(V)$ is a connected graded Hopf algebra for which
the canonical inclusion of $V$ is primitive. For $n\geq2$ and $1\leq i<n$, let
$
 c_i=\id_V^{\otimes(i-1)}\otimes c_{V,V}
             \otimes\id_V^{\otimes(n-i-1)}
 \in\operatorname{End}_{\D}(V^{\otimes n}).
$
Write $s_i=(i,i+1)\in\mathfrak S_n$. For a reduced expression
$\sigma=s_{i_1}\cdots s_{i_r}$, define
\[
 c_\sigma=c_{i_1}\circ\cdots\circ c_{i_r}.
\]
The braid relations and \cite[Theorem~1.7.4]{rootsys} imply that
$c_\sigma$ is independent of the reduced expression.
For the identity permutation, we set
$c_{\mathrm{id}}=\id_{V^{\otimes n}}$, also for $n=0,1$.

Set
$
 T_{\mathrm{sh}}(V)=\bigoplus_{n\geq0}V^{\otimes n},
$
and let $\iota_n:V^{\otimes n}\to T_{\mathrm{sh}}(V)$ be the
canonical inclusion. Define the coproduct by
\[
 \Delta\iota_n
 =\sum_{p+q=n}(\iota_p\otimes\iota_q),
\]
where each summand uses the canonical identification
$V^{\otimes n}\cong V^{\otimes p}\otimes V^{\otimes q}$. Its counit is the projection
onto $V^{\otimes0}=\one$.
The quantum shuffle product is defined by
\[
 m_{T_{\mathrm{sh}}(V)}(\iota_p\otimes\iota_q)
 =\iota_{p+q}\circ
   \left(\sum_{\sigma\in\operatorname{Sh}(p,q)}c_\sigma\right),
 \qquad p,q\geq0,
\]
where $\operatorname{Sh}(p,q)\subseteq\mathfrak S_{p+q}$
consists of the permutations satisfying
\[
 \sigma(1)<\cdots<\sigma(p),
 \qquad
 \sigma(p+1)<\cdots<\sigma(p+q).
\]
The unit is $\iota_0$. These operations make $T_{\mathrm{sh}}(V)$
a connected graded Hopf algebra, called the quantum shuffle algebra,
see \cite[Section~6.5, p.~265]{rootsys}.
The bialgebra identities follow from identities between finite sums
of braid operators, so the construction applies in $\D$.

The universal property of $T(V)$ gives a unique graded algebra morphism
\[
 \mathcal S:T(V)\longrightarrow T_{\mathrm{sh}}(V)
\]
whose restriction to $V$ is $\id_V$. Compatibility with coproducts
and counits follows by comparing algebra morphisms on the generators
$V$, which are primitive in both Hopf algebras, see \
\cite[Proposition~6.4.7]{rootsys}.
Thus $\mathcal S$ is a Hopf algebra morphism.
Iterating the shuffle product gives
\[
 \mathcal S_n=\sum_{\sigma\in\mathfrak S_n}c_\sigma
 \quad(n\geq2),\qquad
 \mathcal S_0=\id_{\one},\quad \mathcal S_1=\id_V,
\]
by the shuffle decomposition in \cite[Corollary~1.8.4]{rootsys}.
The endomorphism $\mathcal S_n$ is called the quantum symmetrizer.

\begin{definition}\label{pre:nichols-definition}
The Nichols algebra of $V$ is the graded Hopf algebra
\[
 \B(V)=\im\mathcal S.
\]
 {Thus $\mathcal S$ factors as}
$
 T(V)\twoheadrightarrow\B(V)\hookrightarrow T_{\mathrm{sh}}(V),
$
and its degree-$n$ component identifies
$\B(V)(n)$ with $\im\mathcal S_n\subseteq V^{\otimes n}$.
\end{definition}

For a morphism $f:V\to W$, naturality of the braiding gives
\[
 f^{\otimes n}\mathcal S_n(V)
 =\mathcal S_n(W)f^{\otimes n}.
\]
The induced morphism
$T_{\mathrm{sh}}(f):T_{\mathrm{sh}}(V)\to T_{\mathrm{sh}}(W)$
therefore restricts to a graded Hopf algebra morphism
\[
 \B(f):\B(V)\longrightarrow\B(W)
\]
whose degree-one component is $f$.
Since $\B(V)$ is generated by its degree-one component,
$\B(f)$ is the unique graded Hopf algebra morphism extending $f$.
For composable morphisms $f:V\to W$ and $g:W\to Z$,
this uniqueness gives
\[
 \B(gf)=\B(g)\B(f),
 \qquad
 \B(\id_V)=\id_{\B(V)}.
\]

\begin{lemma}\label{pre:nichols-properties}
Let $V\in\D$.
\begin{enumerate}
 \item The Hopf algebra $\B(V)$ is connected, is generated by
 $\B(V)(1)=V$, and satisfies $P(\B(V))=V$.
 \item The kernel of $T(V)\longrightarrow\B(V)$ is the largest
 homogeneous Hopf ideal of $T(V)$ with zero components in degrees
 zero and one.
 \item If $V\in\C$, then $\B(V)$ is locally finite.
\end{enumerate}
\end{lemma}

\begin{proof}

Connectedness and generation follow from the definition.
For $n\geq2$, consider the composite
\[
 \B(V)(n)
 \xrightarrow{\Delta_{1,n-1}}
 V\otimes\B(V)(n-1)
 \hookrightarrow V\otimes V^{\otimes(n-1)}
 \cong V^{\otimes n}.
\]
Since $\B(V)$ is a graded subcoalgebra of $T_{\mathrm{sh}}(V)$,
this composite is the canonical inclusion into $V^{\otimes n}$.
Hence $\Delta_{1,n-1}$ is a monomorphism.
Its restriction to $P(\B(V))(n)$ is zero, since a primitive
subobject in degree $n$ has coproduct components only in
bidegrees $(0,n)$ and $(n,0)$.
Thus $P(\B(V))(n)=0$ for $n\geq2$.
Together with connectedness and $\B(V)(1)=V$, this proves~(1).

Let $J$ be a homogeneous Hopf ideal as in (2), and set $Q=T(V)/J$.
The degree-$(1,\ldots,1)$ component of the iterated coproduct
of $T_{\mathrm{sh}}(V)$ is the canonical identification
$V^{\otimes n}\cong V\otimes\cdots\otimes V$.
Since $\mathcal S$ is a coalgebra morphism and is the identity
in degree one, the corresponding component for $T(V)$ is
$\mathcal S_n$.
For $n\geq2$, the quotient morphism and $Q(1)=V$ therefore give
a factorization
\[
 T(V)(n)\longrightarrow Q(n)\longrightarrow V^{\otimes n}
\]
of $\mathcal S_n$. Hence $J(n)\subseteq\ker\mathcal S_n$ for every
$n$. Since $\ker\mathcal S$ is itself a homogeneous Hopf ideal with
zero components in degrees zero and one, this proves (2).

If $V\in\C$, then $V^{\otimes n}\in\C$, and its endomorphism
$\mathcal S_n$ has image in $\C$. This proves (3).
\end{proof}

\begin{lemma}\label{pre:nichols-characterization}
Let $H$ be a connected nonnegatively graded Hopf algebra in $\D$, and
let $V=H(1)$. If $H$ is generated by $V$ and $P(H)=V$, then the identity
of $V$ extends to a graded Hopf algebra isomorphism
\[
 H\cong\B(V).
\]
\end{lemma}

\begin{proof}
By Lemma~\ref{pre:nichols-properties}(2), the canonical epimorphism
$T(V)\to H$ induces a graded Hopf algebra epimorphism
\[
 f:H\longrightarrow\B(V)
\]
which is the identity in degrees zero and one.
Suppose that $J=\ker f$ is nonzero, and choose the least $n\geq2$
such that $J(n)\neq0$. Let $j:J(n)\hookrightarrow H$ be the inclusion.
Since $f$ is a coalgebra morphism and $fj=0$, we have
\[
 (f\otimes f)\Delta_Hj=\Delta_{\B(V)}fj=0.
\]
For $0<i<n$, the restrictions of $f$ to $H(i)$ and $H(n-i)$
are isomorphisms by the choice of $n$. Taking the component of
bidegree $(i,n-i)$ in this equality therefore gives
$\Delta_{i,n-i}j=0$.

Thus $\Delta_Hj$ has only components of bidegrees $(0,n)$ and
$(n,0)$. Since $H$ is connected, $\eta_H:\one\to H(0)$ is an
isomorphism with {inverse} $\varepsilon_H|_{H(0)}$.
Applying the two counit identities to these components gives
\[
 \Delta_Hj=\eta_H\otimes j+j\otimes\eta_H.
\]
Thus $j$ factors through $P(H)(n)=0$, a contradiction.
It follows that $J=0$, so $f$ is an isomorphism.
\end{proof}

\begin{lemma}\label{pre:nichols-antipode}
The antipode of $\B(V)$ is bijective for every $V\in\D$.
\end{lemma}

\begin{proof}
Fix $n\geq2$, and let $w_0\in\mathfrak S_n$ be the longest
permutation. Since $V$ is primitive in $T(V)$, its antipode
restricts to $-\id_V$ on $V$.
Then we have 
$
 S_{T(V)}|_{V^{\otimes n}}
 =-c_{V^{\otimes(n-1)},V}
   \bigl(S_{T(V)}|_{V^{\otimes(n-1)}}\otimes\id_V\bigr).
$
Using $c_{V^{\otimes(n-1)},V}=c_1\cdots c_{n-1}$ and the
reduced expression
$
 w_0=(s_1\cdots s_{n-1})(s_1\cdots s_{n-2})\cdots s_1,
$
induction on $n$ yields
\[
 S_{T(V)}|_{V^{\otimes n}}=(-1)^n c_{w_0}.
\]
We show that $c_{w_0}$ commutes with $\mathcal S_n$.
For $1\leq i<n$, we have
$w_0s_i=s_{n-i}w_0$, and this permutation has length
$\ell(w_0)-1$. The corresponding reduced expressions give
$
 c_{w_0}
 =c_{w_0s_i}c_i
 =c_{n-i}c_{w_0s_i}.
$
It follows that $c_{w_0}c_i=c_{n-i}c_{w_0}$.
Applying this identity to a reduced expression of
$\sigma\in\mathfrak S_n$, we obtain
\[
 c_{w_0}c_\sigma c_{w_0}^{-1}=c_{w_0\sigma w_0}.
\]
Since conjugation by $w_0$ permutes $\mathfrak S_n$,
\[
 c_{w_0}\mathcal S_n c_{w_0}^{-1}
 =\sum_{\sigma\in\mathfrak S_n}c_{w_0\sigma w_0}
 =\mathcal S_n.
\]
Thus both $c_{w_0}$ and $c_{w_0}^{-1}$ preserve
$\ker\mathcal S_n$. Consequently, $(-1)^n c_{w_0}$ induces
an automorphism of
\[
 \B(V)(n)\cong V^{\otimes n}/\ker\mathcal S_n,
\]
whose inverse is induced by $(-1)^n c_{w_0}^{-1}$.
The canonical Hopf algebra projection $T(V)\to\B(V)$ commutes
with antipodes, so this induced automorphism is precisely
$S_{\B(V)}|_{\B(V)(n)}$.
In degrees zero and one, the antipode is $\id_{\one}$ and
$-\id_V$, respectively. Hence it is bijective in every degree,
and therefore bijective on $\B(V)$.
\end{proof}

\begin{lemma}\label{pre:nichols-functoriality}
Let $F$ be an exact braided monoidal functor between braided
abelian categories with exact direct sums.
Assume that their tensor products are exact and preserve
direct sums in each variable, and that $F$ preserves direct sums.
Then there is a natural isomorphism of graded Hopf algebras
\[
 F(\B(V))\cong\B(F(V)).
\]
\end{lemma}

\begin{proof}
The monoidal structure of $F$ and preservation of direct sums
give natural isomorphisms of graded Hopf algebras
$
 F(T(V))\cong T(F(V)),
 \
 F(T_{\mathrm{sh}}(V))\cong T_{\mathrm{sh}}(F(V)).
$
Since $F$ preserves braidings and finite sums, these
isomorphisms identify $F(\mathcal S_V)$ with
$\mathcal S_{F(V)}$.
Exactness of $F$ then gives
\[
 F(\B(V))
 =F(\im\mathcal S_V)
 \cong\im F(\mathcal S_V)
 \cong\im\mathcal S_{F(V)}
 =\B(F(V)).
\]
The induced isomorphism respects the Hopf algebra structures.
For every morphism $f:V\to W$, the tensor and shuffle
morphisms commute with the above identifications.
Passing to images shows that the induced isomorphisms commute
with $F(\B(f))$ and $\B(F(f))$, proving naturality.
\end{proof}

\subsection{Yetter--Drinfeld modules}\label{pre:yet-drinfeld}

Let $A$ be a Hopf algebra in $\D$, with multiplication $\mu_A$
and bijective antipode.
A left Yetter--Drinfeld module over $A$ is an object $X\in\D$
equipped with a left $A$-action $\rho_X:A\otimes X\longrightarrow X$ and a left $A$-coaction $\delta_X:X\longrightarrow A\otimes X$
satisfying
\begin{equation}\label{pre:yd-condition}
\begin{aligned}
 (\mu_A\otimes\rho_X)
  & (\id_A\otimes c_{A,A}\otimes\id_X)
   (\Delta_A\otimes\delta_X)\\
 &=
   (\mu_A\otimes\id_X)
   (\id_A\otimes c_{X,A})
   (\delta_X\rho_X\otimes\id_A)
   (\id_A\otimes c_{A,X})
   (\Delta_A\otimes\id_X).
\end{aligned}
\end{equation}
This is the left Yetter--Drinfeld condition of
\cite[Definition~3.4.1]{rootsys}.
We write ${}_A^A\YD(\D)$ for the category whose objects are
these triples $(X,\rho_X,\delta_X)$ and whose morphisms are
$A$-linear and $A$-colinear morphisms in $\D$.

For $X,Y\in{}_A^A\YD(\D)$, the tensor product $X\otimes Y$
has action and coaction
\begin{align*}
 \rho_{X\otimes Y}
 &=(\rho_X\otimes\rho_Y)
   (\id_A\otimes c_{A,X}\otimes\id_Y)
   (\Delta_A\otimes\id_X\otimes\id_Y),\\
 \delta_{X\otimes Y}
 &=(\mu_A\otimes\id_X\otimes\id_Y)
   (\id_A\otimes c_{X,A}\otimes\id_Y)
   (\delta_X\otimes\delta_Y).
\end{align*}
The unit object is $\one$, with action $\varepsilon_A$ and
coaction $\eta_A$, using the unit constraints.
These structures make ${}_A^A\YD(\D)$ a braided monoidal
category, with braiding
\begin{equation}\label{pre:yd-braiding}
 c^{\YD}_{X,Y}
 =(\rho_Y\otimes\id_X)
   (\id_A\otimes c_{X,Y})
   (\delta_X\otimes\id_Y)
 :X\otimes Y\longrightarrow Y\otimes X,
\end{equation}
see \cite[Proposition~3.2]{Simon}.
The bijectivity of the antipode ensures that
$c^{\YD}_{X,Y}$ is an isomorphism.

We will also use the following characterization of
\eqref{pre:yd-condition}.
Suppose that $X$ is a left $A$-module and a left $A$-comodule.
For every left $A$-module $Y$, formula
\eqref{pre:yd-braiding} defines a morphism
$X\otimes Y\to Y\otimes X$.
Then $X$ is a Yetter--Drinfeld module if and only if all these
morphisms are $A$-linear, where both tensor products carry the
diagonal action above.
It suffices to check this for $Y=A$ with action $\mu_A$,
see \cite[Proposition~3.4.5]{rootsys}.

\begin{lemma}\label{pre:yd-abelian}
The braided monoidal category ${}_A^A\YD(\D)$ is abelian.
The forgetful functor
\[
 {}_A^A\YD(\D)\longrightarrow\D
\]
creates kernels, cokernels, and colimits.
Its tensor product is exact and preserves colimits
in each variable.
\end{lemma}

\begin{proof}

The category ${}_A^A\YD(\D)$ is additive.
Exactness of $A\otimes-$ gives the kernel and cokernel
of any Yetter--Drinfeld morphism unique induced actions
and coactions.
Their defining identities and universal properties follow
from those in $\D$, since tensor products preserve
monomorphisms and epimorphisms.
Thus the forgetful functor creates kernels and cokernels.
Moreover, the inverse of an invertible underlying morphism
is again linear and colinear.
Since the canonical morphism from coimage to image is
an isomorphism in $\D$, it is also an isomorphism in
${}_A^A\YD(\D)$.
Hence this category is abelian.

For a diagram $(X_t)_t$, let $X=\varinjlim_tX_t$ in $\D$,
with structure morphisms $u_t:X_t\to X$.
Preservation of colimits by $A\otimes-$ gives unique
morphisms $\rho_X:A\otimes X\to X$ and
$\delta_X:X\to A\otimes X$ satisfying
\[
 \rho_X(\id_A\otimes u_t)=u_t\rho_{X_t},
 \qquad
 \delta_Xu_t=(\id_A\otimes u_t)\delta_{X_t}.
\]
For $r=0,1,2$, morphisms from $A^{\otimes r}\otimes X$
are determined by their composites with
$\id_{A^{\otimes r}}\otimes u_t$.
The displayed identities therefore reduce the module,
comodule, and Yetter--Drinfeld axioms to those on $X_t$.
They also show that the morphism induced by any compatible
family of Yetter--Drinfeld morphisms is linear and colinear.
Thus the forgetful functor creates colimits.

The tensor product is comleted in $\D$.
Its exactness and preservation of colimits now follow
from the corresponding properties in $\D$.
\end{proof}

\section{ {Dual pairs} and rational Yetter--Drinfeld modules}\label{sec:duality}
We use duality in each degree to construct a braided monoidal
equivalence between rational Yetter--Drinfeld categories.
This equivalence preserves Nichols algebras.

\subsection{Dual pairs}
We continue to work in $\D=\Ind(\C)$, with left actions and left
coactions on Yetter--Drinfeld modules.
For a Hopf algebra $A$ in a braided category $\mathcal C$, we write
${}_A\mathcal C$ and ${}^A\mathcal C$ for its categories of left
modules and left comodules, respectively.
For a braided category $\mathcal C$,
we write $\overline{\mathcal C}$ for the same monoidal category with
braiding $\overline c_{X,Y}=c_{Y,X}^{-1}$. We use the conventions for
Hopf pairings  in \cite[Section~2]{Simon}.
For a bialgebra $R$ in $\mathcal C$, write $R^{\mathrm{cop}}$
for the bialgebra in $\overline{\mathcal C}$ with the same
multiplication and comultiplication $c_{R,R}^{-1}\Delta_R$.
If $R$ has a bijective antipode $S_R$, then $R^{\mathrm{cop}}$
is a Hopf algebra with antipode $S_R^{-1}$.
A Hopf pairing $e:R\otimes S\to\one$ satisfies
\begin{align}
 e(\mu_R\otimes\id_S)
 &=e(\id_R\otimes e\otimes\id_S)
       (\id_R\otimes\id_R\otimes\Delta_S),\label{dual:pairing-first}\\
 e(\id_R\otimes\mu_S)
 &=e(\id_R\otimes e\otimes\id_S)
       (\Delta_R\otimes\id_S\otimes\id_S),\label{dual:pairing-second}\\
 e(\eta_R\otimes\id_S)&=\varepsilon_S,
 \qquad e(\id_R\otimes\eta_S)=\varepsilon_R.\nonumber
\end{align}
Let $R=\bigoplus_{n\geq0}R_n$ and
$S=\bigoplus_{n\geq0}S_n$ be graded Hopf algebras in $\D$,
and let $e:R\otimes S\to\one$ be a Hopf pairing.
We call $(R,S,e)$ a \textbf{dual pair} if
$
 e|_{R_n\otimes S_m}=0 \ (n\ne m),
$
and, for every $n\geq0$, the restriction
$
 e_n:R_n\otimes S_n\longrightarrow\one
$
exhibits $S_n$ as a right dual of $R_n$.

\Needspace{9\baselineskip}
\begin{lemma}\label{dual:pairing}
Let $U\in\C$, and let $U^*$ be a left dual of $U$. Set
$A=\B(U)$ and $B=\B(U^*)$. The evaluation
$U^*\otimes U\to\one$ extends to a dual pair
\[
 p:B\otimes A\longrightarrow\one.
\]
\end{lemma}

\begin{proof}

By Lemma~\ref{pre:nichols-properties}(3), we have
$A_n,B_n\in\C$ for every $n\geq0$.
Let
$
 \operatorname{ev}_n:
 (U^*)^{\otimes n}\otimes U^{\otimes n}\longrightarrow\one
$
be the evaluation obtained from the usual identification
$(U^{\otimes n})^*\cong(U^*)^{\otimes n}$.
The construction in
\cite[Section~5.5 and Remark~5.6]{nicholscatdef}
gives a Hopf pairing
\[
 p^T:T(U^*)\otimes T(U)\longrightarrow\one
\]
extending the evaluation on $U^*\otimes U$.
The pairing vanishes on components of distinct degrees.
By \cite[Section~5.5]{nicholscatdef}, its degree-$(n,n)$
component satisfies
\[
 p_n^T
 =\operatorname{ev}_n
   (\mathcal S_n(U^*)\otimes\id_{U^{\otimes n}})
 =\operatorname{ev}_n
   (\id_{(U^*)^{\otimes n}}\otimes\mathcal S_n(U)),
 \qquad n\geq0,
\]
where $\operatorname{ev}_0=\id_{\one}$.
By uniqueness of transposes with respect to $\operatorname{ev}_n$,
this identity gives
$\mathcal S_n(U)^*=\mathcal S_n(U^*)$
under $(U^{\otimes n})^*\cong(U^*)^{\otimes n}$.
Write the factorization of $\mathcal S_n(U)$ as
\[
 \mathcal S_n(U)=j_nq_n,\qquad
 U^{\otimes n}\xrightarrow{\ q_n\ }A_n
                  \xrightarrow{\ j_n\ }U^{\otimes n}.
\]
Since duality is exact, $j_n^*$ is an epimorphism,
$q_n^*$ is a monomorphism, and
$
 \mathcal S_n(U^*)=\mathcal S_n(U)^*=q_n^*j_n^*.
$
It follows that
$
 \ker j_n^*=\ker\mathcal S_n(U^*).
$
Thus $j_n^*$ induces an isomorphism $B_n\cong A_n^*$.
Naturality of evaluation gives
\[
 \operatorname{ev}_{A_n}(j_n^*\otimes q_n)
 =\operatorname{ev}_n
       (\id_{(U^*)^{\otimes n}}\otimes j_nq_n)
 =p_n^T.
\]
Consequently, $p_n^T$ descends to a pairing
$
 p_n:B_n\otimes A_n\longrightarrow\one
$
which, under the isomorphism $B_n\cong A_n^*$,
is the evaluation $A_n^*\otimes A_n\to\one$.
In particular, $p_n$ exhibits $A_n$ as a right dual of $B_n$.

The pairings $p_n$, together with the zero morphisms in distinct
degrees, define $p:B\otimes A\to\one$.
The Hopf pairing identities descend from those for $p^T$,
since tensor products preserve epimorphisms.
Finally, $\mathcal S_1(U)=\id_U$, so the restriction of $p$
to $U^*\otimes U$ is the given evaluation.
\end{proof}

\subsection{Rational actions and transposition}

Let $R=\bigoplus_{n\geq0}R_n$ be a graded Hopf algebra in $\D$, with
$R_n\in\C$. For each $n\geq0$, choose a right dual ${^*R_n}$ of $R_n$,
with evaluation and coevaluation
$
 \operatorname{ev}_{R_n}:R_n\otimes {^*R_n}\to\one,
 \
 \operatorname{coev}_{R_n}:\one\to {^*R_n}\otimes R_n.
$
Let $\rho:R\otimes X\to X$ be an action.
Under the canonical isomorphism
\[
 \operatorname{Hom}_{\D}(R_n\otimes X,X)
 \cong
 \operatorname{Hom}_{\D}(X,{^*R_n}\otimes X),
\]
the restriction $\rho|_{R_n\otimes X}$ corresponds to the morphism
\[
 t_n
 =(\id_{{^*R_n}}\otimes\rho|_{R_n\otimes X})
   (\operatorname{coev}_{R_n}\otimes\id_X):
 X\longrightarrow {^*R_n}\otimes X.
\]
Define the subobjects
\begin{equation}\label{dual:rational-filtration}
 F_NX=\bigcap_{n>N}\ker t_n\qquad(N\geq0).
\end{equation}
Equivalently, $F_NX$ is the largest subobject of $X$ annihilated by
$R_{>N}=\bigoplus_{n>N}R_n$. We call $X$ \textbf{rational} if the
canonical morphism $\varinjlim_N F_NX\to X$ is an isomorphism.
This definition is independent of the chosen right duals. A
Yetter--Drinfeld module is rational if its underlying module is rational.
We denote this category by $_R^R\YD(\D)_{\rat}$.
\begin{remark}\upshape
In the vector-space case,
\[
 F_NX=\{x\in X\mid R_nx=0\text{ for all }n>N\},
 \qquad
 \varinjlim_NF_NX=\bigcup_{N\geq0}F_NX.
\]
Hence $X$ is rational if and only if, for every $x\in X$,
there exists $N\geq0$ such that $R_nx=0$ for all $n>N$.
This agrees with \cite[Definition~12.2.3]{rootsys}.
\end{remark}
\begin{lemma}\label{dual:rational-closure}
Rational $R$-modules are closed under subobjects, quotient objects,
direct sums and tensor products. The same assertions hold for rational
Yetter--Drinfeld modules. Consequently, ${}_R^R\YD(\D)_{\rat}$ is an
abelian monoidal full subcategory of ${}_R^R\YD(\D)$, and its inclusion
is exact.
\end{lemma}

\begin{proof}

Let $X$ be a rational $R$-module, and let
$i:Y\hookrightarrow X$ be the inclusion of a submodule.
Write $t_n^Y$ and $t_n^X$ for the transposes of the actions.
Since $i$ is $R$-linear, naturality of transposition gives
\[
 t_n^Xi=(\id_{{^*R_n}}\otimes i)t_n^Y.
\]
Exactness of tensor products implies that
$\id_{{^*R_n}}\otimes i$ is a monomorphism. Hence
$\ker t_n^Y=i^{-1}(\ker t_n^X)$.
Taking intersections over $n>N$ gives $F_NY=i^{-1}(F_NX)$.
Since filtered colimits are exact,
$\varinjlim_NF_NY\cong Y$, so $Y$ is rational.
For a quotient
$X\twoheadrightarrow Z$, the image of $F_NX$ is contained in $F_NZ$,
and these images exhaust $Z$. Finite sums of subobjects annihilated by
$R_{>N}$ have the same property. Applying this observation to finite
subsets of a direct sum proves closure under arbitrary direct sums.

Let $X,Y$ be rational $R$-modules. If $X_0\subseteq X$ and
$Y_0\subseteq Y$ are annihilated by $R_{>N}$ and $R_{>M}$,
respectively, then the diagonal action vanishes on
$R_{>N+M}\otimes X_0\otimes Y_0$.
Since tensor products preserve filtered colimits, taking
$X_0=F_NX$ and $Y_0=F_MY$ proves that $X\otimes Y$ is rational.
All arguments apply to the underlying actions of
Yetter--Drinfeld modules. The first two assertions
therefore also prove the last statement.
\end{proof}

\begin{lemma}\label{dual:transposition}
Let $e:B\otimes A\to\one$ be a dual pair of graded
Hopf algebras in $\D$, with homogeneous components in $\C$.
We have a strict monoidal isomorphism
\[
{}^A\D\longrightarrow
 {}_{B^{\mathrm{cop}}}\overline{\D}_{\rat},
 \qquad
 (X,\delta)\longmapsto
 \bigl(X,(e\otimes\id_X)(\id_B\otimes\delta)\bigr).
\]
It acts as the identity on underlying objects and morphisms.
\end{lemma}

\begin{proof}

Let $b_n:\one\to A_n\otimes B_n$ be the coevaluation for $e_n$,
and let $\iota_n:A_n\to A$ and $\pi_n:A\to A_n$ be the
canonical inclusion and projection.
First take a morphism $\rho:B\otimes X\to X$ and define
$F_NX$ by \eqref{dual:rational-filtration}. Suppose that the canonical morphism
$\varinjlim_NF_NX\to X$ is an isomorphism. Write
\[
 \rho_n=\rho|_{B_n\otimes X},
 \qquad
 t_n=(\id_{A_n}\otimes\rho_n)(b_n\otimes\id_X).
\]
By \eqref{dual:rational-filtration}, $t_n$ vanishes on $F_NX$ for
$n>N$. For each $N\geq0$, define
\begin{equation}\label{dual:finite-transpose}
 \delta_N=
 \sum_{n=0}^{N}(\iota_n\otimes\id_X)t_n|_{F_NX}
 :F_NX\longrightarrow A\otimes X.
\end{equation}
By the definition of $F_NX$, we have
$t_{N+1}|_{F_NX}=0$. Hence
\[
 \delta_{N+1}|_{F_NX}
 =\delta_N+
   (\iota_{N+1}\otimes\id_X)t_{N+1}|_{F_NX}
 =\delta_N.
\]
Since the canonical morphism $\varinjlim_NF_NX\to X$ is an
isomorphism, the universal property of the colimit gives
a unique morphism
\[
 \delta:X\longrightarrow A\otimes X
\]
such that $\delta|_{F_NX}=\delta_N$ for every $N\geq0$.
For $N\geq n$, \eqref{dual:finite-transpose} gives
$(\pi_n\otimes\id_X)\delta_N=t_n|_{F_NX}$.
Passing to the colimit over $N\geq n$, we obtain
$(\pi_n\otimes\id_X)\delta=t_n$.

Conversely, let $\delta:X\to A\otimes X$ be a morphism.
For each $N\geq0$, let
\[
 G_NX=
 \ker\bigl(
 X\xrightarrow{\delta}A\otimes X
 \longrightarrow (A/A_{\leq N})\otimes X
 \bigr).
\]
Thus $\delta|_{G_NX}$ factors through $A_{\leq N}\otimes X$.
Since
$\varinjlim_N(A_{\leq N}\otimes X)\cong A\otimes X$,
exactness of filtered colimits gives
$
 \varinjlim_N\bigl((A/A_{\leq N})\otimes X\bigr)=0.
$
Applying the same exactness to the kernels defining $G_NX$,
we obtain
\[
 \varinjlim_NG_NX
 \cong
 \ker\Bigl(
 X\longrightarrow
 \varinjlim_N\bigl((A/A_{\leq N})\otimes X\bigr)
 \Bigr)
 =X.
\]
Define
$
 \rho=(e\otimes\id_X)(\id_B\otimes\delta):
 B\otimes X\longrightarrow X.
$
The pairing $e$ vanishes on $B_{>N}\otimes A_{\leq N}$.
Since $\delta|_{G_NX}$ factors through $A_{\leq N}\otimes X$,
the restriction of $\rho$ to $B_{>N}\otimes G_NX$ is zero.
Consequently, the subobjects defined from $\rho$ by
\eqref{dual:rational-filtration} satisfy $G_NX\subseteq F_NX$.
It follows that $\varinjlim_NF_NX\cong X$.

We check that these constructions are inverse.
If $\rho$ is obtained from $\delta$, the triangle identities give
$
 t_n
 =(\id_{A_n}\otimes\rho_n)(b_n\otimes\id_X)
 =(\pi_n\otimes\id_X)\delta.
$
The reconstructed morphism has the same components $t_n$.
Lemma~\ref{pre:finite-presentability}(2) therefore identifies
it with $\delta$.
Conversely, let $\delta$ be constructed from $\rho$.
Then $(\pi_n\otimes\id_X)\delta=t_n$ for every $n\geq0$.
The other triangle identity gives
$
 (e_n\otimes\id_X)(\id_{B_n}\otimes t_n)=\rho_n.
$
This recovers $\rho$ on each
$B_n\otimes X$, and hence on $B\otimes X$.

The unit and counit conditions correspond, since
\[
 \rho(\eta_B\otimes\id_X)
 =(\varepsilon_A\otimes\id_X)\delta.
\]
We next check associativity and coassociativity.
Fix $r,s\geq0$, and denote the multiplication component by
$\mu_{s,r}:B_s\otimes B_r\to B_{r+s}$.
By \eqref{dual:pairing-first}, its transpose is the component
$A_{r+s}\to A_r\otimes A_s$ of $\Delta_A$.
Thus the transpose of
$\rho_{r+s}(\mu_{s,r}\otimes\id_X)$ is
\[
 (\pi_r\otimes\pi_s\otimes\id_X)
 (\Delta_A\otimes\id_X)\delta.
\]
For the right-hand side, first transpose
$\rho_s(\id_{B_s}\otimes\rho_r)$ along $B_s$.
By naturality of transposition, this gives $t_s\rho_r$.
Transposing this morphism along $B_r$ then gives
\[
 (\id_{A_r}\otimes t_s)t_r
 =
 (\pi_r\otimes\pi_s\otimes\id_X)
 (\id_A\otimes\delta)\delta.
\]
Hence the two sides of associativity correspond to the
degree-$(r,s)$ components of the two sides of coassociativity.
By Lemma~\ref{pre:finite-presentability}(2), two morphisms
into $A\otimes A\otimes X$ are equal if their compositions
with the projection onto $A_r\otimes A_s\otimes X$
agree for every $r,s\geq0$.
Thus $\rho$ is associative if and only if $\delta$ is coassociative.

Let $X$ and $Y$ be left $A$-comodules.
Write $\delta_{X\otimes Y}$ for their tensor product coaction.
By \eqref{dual:pairing-second} and naturality of the braiding,
\begin{align*}
 (e\otimes\id_X\otimes\id_Y)
   (\id_B\otimes\delta_{X\otimes Y})
 =
 (\rho_X\otimes\rho_Y)
 (\id_B\otimes c_{X,B}^{-1}\otimes\id_Y)
\bigl((c_{B,B}^{-1}\Delta_B)\otimes\id_X\otimes\id_Y\bigr).
\end{align*}
The coproduct of $B^{\mathrm{cop}}$ is
$c_{B,B}^{-1}\Delta_B$.
The braiding in $\overline\D$ satisfies
$\overline c_{B,X}=c_{X,B}^{-1}$.
Thus the displayed morphism is the tensor product action
of $B^{\mathrm{cop}}$ in $\overline\D$.
The coaction $\eta_A$ on $\one$ corresponds to the action
$\varepsilon_B$, so the unit object is also preserved.
For a morphism $f:X\to Y$, naturality of transposition gives
\[
 (\id_A\otimes f)\delta_X=\delta_Yf
 \quad\Longleftrightarrow\quad
 f\rho_X=\rho_Y(\id_B\otimes f).
\]
Hence the constructions identify comodule morphisms with
module morphisms.
Together with the preceding checks, this proves the stated
strict monoidal isomorphism.
\end{proof}

\subsection{The equivalence of rational Yetter--Drinfeld categories}

For the rest of this section, let $p:B\otimes A\to\one$ be a
dual pair of $\mathbb{N}_0$-graded Hopf algebras with $A_n,B_n\in\C$.
Assume that their antipodes are bijective. Write
\begin{align}
 p^+&=p\,c_{A,B}(S_A\otimes S_{B}),\label{dual:pplus}\\
 q&=p^+(\id_A\otimes S_{B}^{-1}).\label{dual:q}
\end{align}
By \cite[Proposition~3.3.8(2),(3)]{rootsys}, $p^+$ is a Hopf
pairing in $\D$, and $q$ is a Hopf pairing of
$A^{\mathrm{cop}}$ and $B^{\mathrm{cop}}$ in $\overline\D$.
Furthermore, 
$p^+$ and $q$ are dual pairs of $\mathbb{N}_0$-graded Hopf algebras.
Lemma~\ref{dual:transposition} gives strict monoidal isomorphisms
\begin{equation}\label{dual:D1-D2}
 D_1:{}^A\D\longrightarrow
 {}_{B^{\mathrm{cop}}}\overline\D_{\rat},
 \qquad
 D_2:{}^{B^{\mathrm{cop}}}\overline\D
       \longrightarrow {}_A\D_{\rat},
\end{equation}
defined using $p$ and $q$, respectively.

\begin{lemma}\label{dual:gamma}
The assignments in \eqref{dual:D1-D2} induce a strict braided monoidal
isomorphism
\[
 \Gamma:\overline{{}_A^A\YD(\D)_{\rat}}\longrightarrow
 {}_{B^{\mathrm{cop}}}^{B^{\mathrm{cop}}}
       \YD(\overline\D)_{\rat}.
\]
For $(X,\rho,\delta)$, the new structures $(\bar\rho,\bar\delta)$ are
determined by
\[
 D_1(X,\delta)=(X,\bar\rho),
 \qquad D_2(X,\bar\delta)=(X,\rho).
\]
\end{lemma}

\begin{proof}
 
Let $V$ have a rational left $A$-action $\rho$ and a left
$A$-coaction $\delta$. Let $\bar V=(V,\bar\rho,\bar\delta)$ be defined by
$D_1$ and $D_2$ as in the statement.
We show that $V$ satisfies the Yetter--Drinfeld condition
if and only if $\bar V$ does.

For a left $A$-module $Z$, let
\[
 b^A_{V,Z}
 =(\rho_Z\otimes\id_V)
   (\id_A\otimes c_{V,Z})(\delta\otimes\id_Z)
 :V\otimes Z\longrightarrow Z\otimes V.
\]
By \cite[Proposition~3.4.5 $(1)\Leftrightarrow(2)$]{rootsys}, $V$ is a
Yetter--Drinfeld module if and only if $b^A_{V,Z}$ is
$A$-linear for every left $A$-module $Z$.
We show that testing all rational $Z$ is also sufficient.

Suppose that $b^A_{V,Z}$ is $A$-linear for every rational
$A$-module $Z$. For $N\geq0$, let
\[
 Z_N=A/A_{>N},
 \qquad q_N:A\longrightarrow Z_N
\]
be the quotient of the regular module.
The grading gives $A_{>N}\cdot Z_N=0$, so $Z_N$ is rational.
Naturality gives
\[
 (q_N\otimes\id_V)b^A_{V,A}
 =b^A_{V,Z_N}(\id_V\otimes q_N).
\]
Both $q_N$ and $b^A_{V,Z_N}$ are $A$-linear.
Thus the composite $(q_N\otimes\id_V)b^A_{V,A}$ is
$A$-linear, which gives
\[
 (q_N\otimes\id_V)b^A_{V,A}\rho_{V\otimes A}
 =
 (q_N\otimes\id_V)\rho_{A\otimes V}
       (\id_A\otimes b^A_{V,A}).
\]
Here both tensor products carry the diagonal action.
For every $n\geq0$, the projection $A\otimes V\to A_n\otimes V$
factors through $q_n\otimes\id_V$, since
$Z_n\cong\bigoplus_{r=0}^n A_r$ in $\D$.
By Lemma~\ref{pre:finite-presentability}(2), these projections
jointly detect equality. Hence
\[
 b^A_{V,A}\rho_{V\otimes A}
 =\rho_{A\otimes V}(\id_A\otimes b^A_{V,A}),
\]
so $b^A_{V,A}$ is $A$-linear.
\cite[Proposition~3.4.5 $(1)\Leftrightarrow(3)$]{rootsys} shows that $V$ is a
Yetter--Drinfeld module.

For a rational $A$-module $Z$, write $\bar Z=D_2^{-1}(Z)$.
Since the antipode of $A$ is bijective, $b^A_{V,Z}$ is
an isomorphism by \cite[Proposition~3.4.8]{rootsys}.
Moreover, \cite[Lemma~3.12]{reflection2} gives
\begin{equation}\label{dual:candidate-braiding}
 (b^A_{V,Z})^{-1}
 =b^{B^{\mathrm{cop}}}_{\bar Z,\bar V}.
\end{equation}
The right-hand side is formed in $\overline\D$.
This identity holds for the actions and coactions defined
by $D_1$ and $D_2$, before imposing the
Yetter--Drinfeld condition.
Its proof uses \eqref{dual:q} and naturality of the braiding.

The morphism $b^A_{V,Z}$ is $A$-linear if and only if its
inverse is $A$-linear.
Since $D_2$ is a strict monoidal isomorphism, this is
equivalent to $B^{\mathrm{cop}}$-colinearity of
$b^{B^{\mathrm{cop}}}_{\bar Z,\bar V}$.
As $Z$ runs through the rational $A$-modules,
$\bar Z$ runs through all left $B^{\mathrm{cop}}$-comodules
in $\overline\D$.
By \cite[Proposition~3.4.5, (1)$\Leftrightarrow$(4)]{rootsys},
these morphisms are colinear for all $\bar Z$ if and only
if $\bar V$ satisfies the Yetter--Drinfeld condition.
Together with the argument above, this proves
the required equivalence.

The inverse correspondences supplied by $D_1$ and $D_2$
therefore restrict to Yetter--Drinfeld modules.
They also identify the corresponding module and comodule
morphisms. Hence $\Gamma$ is an isomorphism of categories.
It is strict monoidal because $D_1$ and $D_2$ are strict
monoidal.

Finally, for $V,W\in{}_A^A\YD(\D)_{\rat}$,
\eqref{dual:candidate-braiding} gives
\[
 b^{B^{\mathrm{cop}}}_{\bar V,\bar W}
 =(b^A_{W,V})^{-1}.
\]
The right-hand side is the braiding from $V\otimes W$
to $W\otimes V$ in $\overline{{}_A^A\YD(\D)_{\rat}}$.
Thus $\Gamma$ is braided.
\end{proof}

\begin{theorem}\label{dual:omega}
Let $(B,A,p)$ be a dual pair of nonnegatively graded Hopf
algebras in $\D$. Assume that $A_n,B_n\in\C$ for every
$n\geq0$ and that both antipodes are bijective.
Then there is an exact braided monoidal equivalence
\[
 \Omega:{}_A^A\YD(\D)_{\rat}
 \longrightarrow {}_B^B\YD(\D)_{\rat}
\]
which is the identity on underlying objects and morphisms.

For $X=(X,\rho,\delta)$, let $\delta_0:X\to B\otimes X$
be the unique morphism satisfying
\begin{equation}\label{dual:delta-zero}
 \rho=(p^+\otimes\id_X)(\id_A\otimes\delta_0),
\end{equation}
where $p^+$ is given by \eqref{dual:pplus}.
The action and coaction on $\Omega(X)$ are
\begin{align}
 \rho^\Omega
   &=(p\otimes\id_X)(\id_B\otimes\delta),
       \label{dual:omega-action}\\
 \delta^\Omega
   &=c_{X,B}c_{B,X}(S_B^2\otimes\id_X)\delta_0.
       \label{dual:omega-coaction}
\end{align}
The monoidal structure is
\begin{equation}\label{dual:tensorator}
 \begin{gathered}
 \beta_{X,Y}=c^{\YD}_{Y,X}c_{Y,X}^{-1}:
 \Omega(X)\otimes\Omega(Y)\longrightarrow\Omega(X\otimes Y),
 \\
 \beta_0=\id_{\one},
 \end{gathered}
\end{equation}
where $c^{\YD}$ is the braiding of ${}_A^A\YD(\D)$
in \eqref{pre:yd-braiding}.
\end{theorem}

\begin{proof}
We construct $\Omega$ from $\Gamma$ by passing through
right-right Yetter--Drinfeld modules over $B$.
Let
$ \YD_B^B(\D)_{\rat}$
be the category whose underlying right $B$-modules are
rational, using the right-handed version of
\eqref{dual:rational-filtration}.
Its braiding is
$
 c^{\mathcal R}_{X,Y}
 =(\id_Y\otimes\rho_X^r)(c_{X,Y}\otimes\id_B)
   (\id_X\otimes\delta_Y^r),
$
{as in \cite[Remark~3.5]{Simon}}.
By \cite[Theorem~3.4.16]{rootsys}, we have the strict braided monoidal equivalence
\[
 F_{lr}:
 {}_{B^{\mathrm{cop}}}^{B^{\mathrm{cop}}}
       \YD(\overline\D)_{\rat}
 \longrightarrow\overline{\YD_B^B(\D)_{\rat}},
\]
given by
$
 F_{lr}(X,\lambda,\nu)
 =\bigl(
 X,\lambda c_{X,B}(\id_X\otimes S_B),c_{B,X}\nu
 \bigr).
$
By  {\cite[Theorem~3.16]{Simon}}, we have the
braided monoidal equivalence
\[
 F_{rl}:\YD_B^B(\D)_{\rat}
 \longrightarrow{}_B^B\YD(\D)_{\rat},
\]
given by
$
 F_{rl}(X,\rho^r,\delta^r)
 =\bigl(
 X,\rho^r c_{X,B}^{-1}(S_B^{-1}\otimes\id_X),
 (S_B\otimes\id_X)c_{X,B}\delta^r
 \bigr).
$
Its tensor constraint is
\[
 \alpha_{X,Y}=c^{\mathcal R}_{Y,X}c_{Y,X}^{-1}:
 F_{rl}(X)\otimes F_{rl}(Y)
 \longrightarrow F_{rl}(X\otimes Y),
\]
and its unit constraint is the identity.
Both functors act as the identity on underlying morphisms.

These functors and their inverses preserve rationality.
Indeed, the corresponding left and right actions are related by
\[
 \rho^r=\lambda c_{X,B}(\id_X\otimes S_B),
 \qquad
 \lambda=\rho^r c_{X,B}^{-1}(S_B^{-1}\otimes\id_X).
\]
The antipode and its inverse preserve each $B_n$.
Hence, for every subobject $X_0\subseteq X$, naturality
of the braiding gives
\[
 \lambda|_{B_{>N}\otimes X_0}=0
 \quad\Longleftrightarrow\quad
 \rho^r|_{X_0\otimes B_{>N}}=0.
\]
Thus $X$ is rational as a left $B$-module if and only if
it is rational as a right $B$-module.
Hence both functors and their inverses restrict to
the rational subcategories.

We compose the three functors on their underlying monoidal
categories and set
\[
 \Omega=F_{rl}F_{lr}\Gamma.
\]
By Lemma~\ref{dual:gamma}, $F_{lr}\Gamma$ is a strict
braided monoidal equivalence
$
 \overline{{}_A^A\YD(\D)_{\rat}}
 \longrightarrow\overline{\YD_B^B(\D)_{\rat}}.
$
It preserves the reversed braidings, and hence also
the original braidings.
Since $F_{rl}$ is braided, $\Omega$ is a braided
monoidal equivalence.  {Since both categories are abelian by
Lemma~\ref{dual:rational-closure}, the equivalence $\Omega$ is exact.}

For $X=(X,\rho,\delta)$, write
$\Gamma(X)=(X,\bar\rho,\bar\delta)$.
By the definitions of $D_1$ and $D_2$,
\[
 \bar\rho=(p\otimes\id_X)(\id_B\otimes\delta),
 \qquad
 \rho=(q\otimes\id_X)(\id_A\otimes\bar\delta).
\]
Lemma~\ref{dual:transposition}, applied to $p^+$, gives the
unique morphism $\delta_0$ in \eqref{dual:delta-zero}.
Since $q=p^+(\id_A\otimes S_B^{-1})$, the second identity becomes
$
 \rho=(p^+\otimes\id_X)
\bigl(\id_A\otimes((S_B^{-1}\otimes\id_X)\bar\delta)\bigr).
$
Uniqueness gives $\bar\delta=(S_B\otimes\id_X)\delta_0$.
The formulas for $F_{lr}$ and $F_{rl}$ and naturality of the
braiding now give
\begin{align*}
 \rho^\Omega
 &=\bar\rho\,c_{X,B}(\id_X\otimes S_B)
       c_{X,B}^{-1}(S_B^{-1}\otimes\id_X)
 =\bar\rho,\\
 \delta^\Omega
 &=(S_B\otimes\id_X)c_{X,B}c_{B,X}\bar\delta
 =c_{X,B}c_{B,X}(S_B^2\otimes\id_X)\delta_0.
\end{align*}
These are \eqref{dual:omega-action}
and \eqref{dual:omega-coaction}.

Since $\Gamma$ and $F_{lr}$ are strict monoidal,
the tensor constraint of $\Omega$ is
$
 \beta_{X,Y}
 =\alpha_{F_{lr}\Gamma(X),F_{lr}\Gamma(Y)}.
$
The composite $F_{lr}\Gamma$ preserves the original
braidings and leaves underlying morphisms unchanged.
Thus
$
 c^{\mathcal R}_{F_{lr}\Gamma(Y),F_{lr}\Gamma(X)}
 =c^{\YD}_{Y,X}.
$
The formula for $\alpha$ now gives
\[
 \beta_{X,Y}=c^{\YD}_{Y,X}c_{Y,X}^{-1},
\]
which is \eqref{dual:tensorator}.
The unit constraint is also the identity.
\end{proof}

For a Hopf algebra $L$ in ${}_A^A\YD(\D)_{\rat}$, the multiplication
and comultiplication transported by $\Omega$ are
\[
m_{\Omega(L)}=\Omega(m_L)\beta_{L,L},
 \qquad
 \Delta_{\Omega(L)}=\beta_{L,L}^{-1}\Omega(\Delta_L).
\]

\begin{corollary}\label{dual:nichols}
If $X\in{}_A^A\YD(\D)_{\rat}$, its Nichols algebra in
${}_A^A\YD(\D)$ is rational, and there is a natural isomorphism of
graded Hopf algebras
\[
 \Omega\bigl(\B(X)\bigr)
   \simeq\B\bigl(\Omega(X)\bigr).
\]
\end{corollary}

\begin{proof}

By Lemma~\ref{dual:rational-closure}, the tensor algebra
of $X$ and its Nichols algebra are rational.
The inclusions of the rational subcategories into the full
Yetter--Drinfeld categories are exact braided monoidal
functors and preserve direct sums.
By Theorem~\ref{dual:omega}, $\Omega$ is exact and braided
monoidal. It also preserves direct sums, since it is an
equivalence.
Applying Lemma~\ref{pre:nichols-functoriality} to $\Omega$
and to the two inclusions gives the stated natural
isomorphism.
\end{proof}

\section{Coinvariants of Nichols algebras}
\label{sec:coinvariants}
We identify the coinvariants associated with a direct-sum
decomposition of the generating object as a Nichols algebra
in a Yetter--Drinfeld category. We also prove that the module
generated by a simple summand of the complement under the
adjoint action is simple in that category.

For a left comodule $(X,\delta_X)$ over a Hopf algebra $A$ in $\D$,
its subobject of left coinvariants is
\[
 \prescript{\co A}{}{X}
 =\ker(\delta_X-\eta_A\otimes\id_X).
\]

\begin{lemma}\label{coinv:nonzero}
Let $A=\bigoplus_{n\geq 0}A_n$ be a connected $\mathbb{N}_0$-graded Hopf algebra
in $\D$. If $X$ is a nonzero left $A$-comodule, then
$\prescript{\co A}{}{X}\neq 0$.
\end{lemma}

\begin{proof}

 Let $A_+=\bigoplus_{n>0}A_n$ and let $q:A\to A_+$ be the projection.
For the coaction $\delta_X:X\to A\otimes X$, set
\[
 \overline\delta_X=(q\otimes\id_X)\delta_X,
 \qquad
 \overline\delta_X^{\,(r+1)}
 =\bigl(\id_{A_+}^{\otimes r}\otimes\overline\delta_X\bigr)
   \overline\delta_X^{\,(r)}
 \quad(r\geq1),
\]
where $\overline\delta_X^{\,(1)}=\overline\delta_X$.
Thus $\overline\delta_X^{\,(r)}$ has target $A_+^{\otimes r}\otimes X$.
 Let $C_nX=\ker\overline\delta_X^{\,(n+1)}$ for $n\geq0$.
The counit identity gives $C_0X=\prescript{\co A}{}{X}$,
and coassociativity gives
\[
 \overline\delta_X^{\,(r)}
 = (q^{\otimes r}\otimes\id_X)
   (\Delta_A^{(r-1)}\otimes\id_X)\delta_X
 \qquad(r\geq 1).
\]
Here $\Delta_A^{(r-1)}:A\to A^{\otimes r}$ is the iterated
coproduct, with $\Delta_A^{(0)}=\id_A$.
Since tensor products preserve direct sums, we have
$
 A\otimes X
 \cong\varinjlim_N
 \left(\bigoplus_{d=0}^{N}A_d\otimes X\right).
$
Let $f:C\to X$ with $C\in\C$.
Finite presentability of $C$ implies that $\delta_Xf$
factors through $\bigoplus_{d=0}^{N}A_d\otimes X$
for some $N$.
For $0\leq d\leq N$, the restriction of $\Delta_A^{(N)}$ to
$A_d$ factors through
\[
 \bigoplus_{\substack{d_1,\ldots,d_{N+1}\geq0\\
                     d_1+\cdots+d_{N+1}=d}}
 A_{d_1}\otimes\cdots\otimes A_{d_{N+1}}.
\]
In each summand, at least one $d_i$ is zero.
Since $q|_{A_0}=0$, applying $q^{\otimes(N+1)}$
annihilates every such summand.
The displayed formula for the iterated reduced coaction
therefore gives
\[
 \overline\delta_X^{\,(N+1)}f=0.
\]
By the definition of $C_NX$, the morphism $f$ factors
through $C_NX\hookrightarrow X$.
The kernels $C_nX$ form an increasing sequence.
Exactness of filtered colimits makes the induced morphism
$\varinjlim_n C_nX\to X$ a monomorphism.
It is an epimorphism because every morphism $C\to X$ with
$C\in\C$ factors through it and the objects of $\C$ generate $\D$.
Hence
\[
 X=\bigcup_{n\geq0}C_nX.
\]
By the definition of the iterates and exactness of $A_+\otimes-$,
\[
 \overline\delta_X(C_nX)\subseteq A_+\otimes C_{n-1}X
 \qquad(n\geq 1).
\]
Since $X\neq0$, choose the least $n$ such that $C_nX\neq0$.
If $n>0$, the last inclusion gives
$C_nX\subseteq\ker\overline\delta_X=C_0X$, a contradiction.
Thus $\prescript{\co A}{}{X}=C_0X\neq0$.
\end{proof}

Let $U,W\in\C$ and let
$
 H=\B(U),\ L=\B(U\oplus W).
$
The projection $U\oplus W\to U$ and the inclusion $U\to U\oplus W$
induce Hopf algebra morphisms
$
 p:L\longrightarrow H,\ i:H\longrightarrow L,
 \ pi=\id_H.
$
Let $K=L^{\co H}$ be the subobject of right coinvariants:
\[
 K=\ker\bigl((\id_L\otimes p)\Delta_L
                   -\id_L\otimes\eta_H\bigr).
\]

Let $j:K\hookrightarrow L$ be the kernel inclusion.
The category $\D$ has equalizers, and $H\otimes-$ preserves them by
exactness. The antipodes of $H$ and $L$ are bijective by
Lemma~\ref{pre:nichols-antipode}. Hence the Radford projection theorem
\cite[Theorem~3.10]{Simon} gives a Hopf algebra structure on $K$ in
${}_H^H\YD(\D)$ and a Hopf algebra isomorphism
\[
 \mu:K\#H\longrightarrow L,
 \qquad \mu=m_L(j\otimes i).
\]
We denote the $H$-action and coaction of $K$ by $\rho_K$ and $\delta_K$.
Recall that the braided adjoint action of $L$ on itself is
\[
 \ad_L
 =m_L(m_L\otimes S_L)
   (\id_L\otimes c_{L,L})(\Delta_L\otimes\id_L)
 :L\otimes L\longrightarrow L.
\]
The $H$-action $\rho_K:H\otimes K\to K$ is the restriction
of this action along $i\otimes j$:
\[
 j\rho_K=\ad_L(i\otimes j).
\]
The Radford projection $r:L\to K$ satisfies
\begin{equation}\label{coinv:projection}
 jr=m_L(\id_L\otimes iS_Hp)\Delta_L,
 \qquad r\mu=\id_K\otimes\varepsilon_H.
\end{equation}
The coalgebra and comodule structures are given by
\begin{equation}\label{coinv:structures}
 \Delta_K=(r\otimes r)\Delta_Lj,
 \qquad
 (\id_H\otimes j)\delta_K=(p\otimes\id_L)\Delta_Lj.
\end{equation}
The primitive subobject of $K$ is
\[
 P(K)=\ker\bigl(\Delta_K-\id_K\otimes\eta_K
                         -\eta_K\otimes\id_K\bigr).
\]
This kernel is taken in ${}_H^H\YD(\D)$.
Since the forgetful functor to $\D$ creates kernels,
its underlying object is also the kernel of the same
morphism in $\D$.
\begin{lemma}\label{coinv:primitive}
Let $v:Z\hookrightarrow K$ be a monomorphism in $\D$
which factors through $P(K)$. Then
\[
\Delta_Ljv
   =jv\otimes\eta_L+(i\otimes j)\delta_Kv.
\]
Consequently, if $Z$ is an $H$-subcomodule of $P(K)$, then
$j(\prescript{\co H}{}{Z})\subseteq P(L)$.
\end{lemma}

\begin{proof}

Since $v$ factors through $P(K)$ and $\eta_K$ is $H$-colinear,
the bosonization coproduct gives
\[
 \begin{aligned}
 \Delta_Ljv
 =(\mu\otimes j)(\id_K\otimes\delta_K)\Delta_Kv
 &=(\mu\otimes j)
   \bigl(v\otimes\eta_H\otimes\eta_K
         +\eta_K\otimes\delta_Kv\bigr)\\
 &=jv\otimes\eta_L+(i\otimes j)\delta_Kv.
 \end{aligned}
\]
The last equality uses
$\mu(\id_K\otimes\eta_H)=j$,
$\mu(\eta_K\otimes\id_H)=i$ and $j\eta_K=\eta_L$.

If $Z$ is an $H$-subcomodule, let
$u:\prescript{\co H}{}{Z}\hookrightarrow Z$ be the inclusion.
Colinearity of $v$ gives
\[
 \delta_Kvu=(\id_H\otimes v)\delta_Zu
           =\eta_H\otimes vu.
\]
Composing the preceding coproduct identity with $u$ yields
\[
 \Delta_Ljvu=jvu\otimes\eta_L+\eta_L\otimes jvu.
\]
Thus $jvu$ factors through $P(L)$.
\end{proof}

The inclusion $W\to L$ factors through $K$, since $W$ is primitive in
$L$ and $p|_W=0$. Define
\[
 Q=(\ad H)(W)
   :=\im\bigl(H\otimes W\hookrightarrow H\otimes K
                    \xrightarrow{\rho_K}K\bigr).
\]

\begin{theorem}\label{coinv:nichols}
The action and coaction of $K$ restrict to $Q$, and there is a unique
isomorphism
\[
 \B(Q)\longrightarrow K
\]
of Hopf algebras in ${}_H^H\YD(\D)$ whose restriction to $Q$ is the
inclusion.
\end{theorem}

\begin{proof}

Give $U$ degree zero and $W$ degree one.
The quantum symmetrizers preserve this grading, so it induces
a Hopf algebra grading
$
 L=\bigoplus_{n\geq0}L[n].
$
Place $H$ in degree zero. The maps $p$ and $i$ are graded.
The morphism defining $K$ is therefore graded, so its kernel has
the grading
\[
 K=\bigoplus_{n\geq0}K[n],
 \qquad K[n]=K\cap L[n],
\]
with $j$ graded. Formula~\eqref{coinv:projection} shows that
$jr$ is graded. Since $j$ is a monomorphism, $r$ is graded as well.
The formulas for the $H$-action and~\eqref{coinv:structures}
now make $K$ a graded Hopf algebra in ${}_H^H\YD(\D)$.
The inclusion $i$ identifies $H$ with $L[0]$, so the unit identifies
\[
 \one\cong H^{\co H}\cong K[0].
\]
Indeed, for any morphism $f:Y\to H$ satisfying
$\Delta_Hf=f\otimes\eta_H$, the counit identity
gives $f=\eta_H\varepsilon_Hf$.
Thus $K$ is connected.

We first prove that $K$ is generated by $Q$.
The $H$-action preserves the grading, so $Q\subseteq K[1]$.
Associativity of the action makes $Q$ an $H$-submodule.
Let $R\subseteq K$ be the image in $\D$ of the graded algebra
morphism $T(Q)\to K$ induced by the inclusion of $Q$.
Since $K$ is an $H$-module algebra, multiplication is
$H$-linear for the diagonal action on tensor powers.
Hence the products of copies of $Q$, and therefore $R$,
are stable under the $H$-action.

The smash-product multiplication on $K\#H$ consequently
restricts to $R\otimes H$.
Its image under $\mu$ is a subalgebra of $L$ containing
$i(H)$ and $j(Q)$.
The unit of $H$ gives $W\subseteq Q$, so this image contains
$U\oplus W$.
Since $L$ is generated by $U\oplus W$, the restriction
\[
 \mu|_{R\otimes H}:R\otimes H\longrightarrow L
\]
is an epimorphism.
Also, \eqref{coinv:projection} gives $rj=\id_K$, so
$r\mu|_{R\otimes H}$ is an epimorphism onto $K$.
By the same identity, this composite factors as
\[
 R\otimes H
 \xrightarrow{\id_R\otimes\varepsilon_H}
 R\hookrightarrow K.
\]
Thus the inclusion $R\hookrightarrow K$ is an epimorphism,
and $R=K$.

Since $R$ is generated by $Q$ in degree one, we obtain
\[
 K[1]=R[1]=Q.
\]
Each homogeneous component of $K$ is a Yetter--Drinfeld
subobject, so this equality gives the required
Yetter--Drinfeld structure on $Q$.
Connectedness of $K$ also gives $Q\subseteq P(K)$.

It remains to prove that $P(K)=Q$.
The primitive subobject $P(K)$ is graded, and its
degree-zero component is zero because $K[0]=\one$.
For $n\geq2$, let
\[
 Z=P(K)\cap K[n].
\]
Both $P(K)$ and $K[n]$ are Yetter--Drinfeld subobjects
of $K$, so $Z$ is an $H$-subcomodule.
Suppose that $Z\neq0$.
With its standard grading, $H=\B(U)$ is connected.
Lemma~\ref{coinv:nonzero} therefore gives
$\prescript{\co H}{}{Z}\neq0$.
Lemma~\ref{coinv:primitive} then implies
\[
 0\neq j\bigl(\prescript{\co H}{}{Z}\bigr)
 \subseteq P(L)\cap L[n].
\]
But $P(L)=U\oplus W$ by
Lemma~\ref{pre:nichols-properties}(1), and these two
summands have degrees zero and one.
Thus $P(L)\cap L[n]=0$, a contradiction.
Consequently,
\[
 P(K)=K[1]=Q.
\]
Lemma~\ref{pre:nichols-characterization}, applied in
${}_H^H\YD(\D)$, gives the required Hopf algebra isomorphism
\[
 \B(Q)\xrightarrow{\sim}K
\]
whose restriction to $Q$ is the inclusion.
It is unique because $\B(Q)$ is generated by $Q$.
\end{proof}

Suppose that $W=\bigoplus_{j\in J}M_j$, where $J$ is finite
and each $M_j$ is simple in $\C$.
For $n\geq0$, let $H(n)$ be the standard degree-$n$
component of $H$.
Define $(\ad U)^n(M_j)$ to be the image of the restricted
action $H(n)\otimes M_j\to K$.
In particular, $(\ad U)^0(M_j)=M_j$.

Since $H$ is generated by $U$, multiplication gives an
epimorphism $U^{\otimes n}\otimes M_j\to H(n)\otimes M_j$.
Associativity of the adjoint action shows that the composite
\begin{equation}\label{coinv:adjoint-factorization}
 U^{\otimes n}\otimes M_j
 \twoheadrightarrow H(n)\otimes M_j
 \longrightarrow K\hookrightarrow L
\end{equation}
is obtained by applying the adjoint action of $U$
successively $n$ times.
Thus the above definition agrees with the iterated
adjoint action.

The grading of $T(U\oplus W)$ by the number of factors
from $U$ and each $M_j$ descends to $L$, since quantum
symmetrizers preserve these numbers.
The subobject $(\ad U)^n(M_j)$ has $U$-degree $n$,
$M_j$-degree $1$, and all other degrees zero.
Hence the sum of these subobjects over $(n,j)$ is direct.
Using $H=\bigoplus_{n\geq0}H(n)$, we obtain
\begin{equation}\label{coinv:adjoint-components}
 Q_j:=(\ad H)(M_j)=\bigoplus_{n\geq0}(\ad U)^n(M_j),
 \qquad
 Q=\bigoplus_{j\in J}Q_j.
\end{equation}

\begin{proposition}\label{coinv:simple}
Each $Q_j$ is simple in ${}_H^H\YD(\D)$, and
\[
 \prescript{\co H}{}{Q_j}=M_j.
\]
\end{proposition}

\begin{proof}
By Theorem~\ref{coinv:nichols} and
\eqref{coinv:adjoint-components}, we have
\[
 P(K)=K[1]=Q=\bigoplus_{j\in J}Q_j.
\]
In the multigrading considered above, $Q_j$ is the part
of $K[1]$ of degree one in $M_j$ and degree zero in
each $M_k$ with $k\neq j$.
The $H$-action and coaction preserve these degrees.
Hence each $Q_j$ is a Yetter--Drinfeld subobject of $P(K)$.

We first identify its left coinvariants.
The unit of $H$ gives $M_j\subseteq Q_j$.
Since $P(L)=U\oplus\bigoplus_{k\in J}M_k$, comparison
of multidegrees gives
$
 j(Q_j)\cap P(L)=M_j
$
as subobjects of $L$.
Lemma~\ref{coinv:primitive} therefore implies
\[
 j\bigl(\prescript{\co H}{}{Q_j}\bigr)
 \subseteq j(Q_j)\cap P(L)=M_j.
\]
Thus $\prescript{\co H}{}{Q_j}\subseteq M_j$.

For the reverse inclusion, let $u_j:M_j\hookrightarrow K$
be the inclusion.
Since $ju_j$ is primitive in $L$ and $pju_j=0$,
\eqref{coinv:structures} gives
\[
 \begin{aligned}
 (\id_H\otimes j)\delta_Ku_j
 =(p\otimes\id_L)\Delta_Lju_j
 =\eta_H\otimes ju_j.
 \end{aligned}
\]
The morphism $\id_H\otimes j$ is a monomorphism,
so $\delta_Ku_j=\eta_H\otimes u_j$.
Since $M_j\subseteq Q_j$, this proves
\[
 \prescript{\co H}{}{Q_j}=M_j.
\]

Let $X\subseteq Q_j$ be a nonzero Yetter--Drinfeld subobject.
Lemma~\ref{coinv:nonzero}, applied to the standard grading of $H$,
and the inclusion $X\hookrightarrow Q_j$ give
\[
 0\neq\prescript{\co H}{}{X}
 \subseteq\prescript{\co H}{}{Q_j}=M_j.
\]
Since $M_j$ is simple in $\D$ by
Lemma~\ref{pre:simple-ind}, the inclusion
$M_j\hookrightarrow Q_j$ factors through $X$.
As $X$ is stable under the $H$-action, the defining morphism
$H\otimes M_j\to Q_j$ also factors through $X$.
Its image is $Q_j$, so $X=Q_j$.
\end{proof}

\section{Reflections of Nichols algebras}\label{sec:reflection}

We prove the reflection theorem and show that it will give rise to the structure of a semi-Cartan graph, extending
\cite{reflection2}, see also \cite{HSdualpair,rootsys}.

\subsection{\texorpdfstring{$\mathbb{Z}$-graded Yetter--Drinfeld modules}
{Z-graded Yetter--Drinfeld modules}}

 {We use} the grading  { convention of Section~\ref{pre:ambient},
allowing degrees in} $  \mathbb Z   $.
Let $A=  \bigoplus_{r\geq0}   A_r$ be a graded Hopf algebra in $\D$,
and let $(X,\rho,\delta)\in{}_A^A\YD(\D)$.
 {A} grading $X=\bigoplus_{d\in\mathbb Z}X_d$  {of its underlying
object is {compatible}} with the Yetter--Drinfeld structure if
\[
 \rho(A_r\otimes X_d)\subseteq X_{r+d},
 \qquad
 \delta(X_d)\subseteq
 \bigoplus_{r\geq0}A_r\otimes X_{d-r}.
\]

\begin{lemma}\label{dual:grading}
Let $(B,A,p)$ and $\Omega$ be as in
Theorem~\ref{dual:omega}.
Suppose that $X\in{}_A^A\YD(\D)_{\rat}$ has a compatible
$\mathbb Z$-grading.
Then $\Omega(X)$ has a compatible grading given by
\[
 \Omega(X)_d=X_{-d}.
\]
For two such objects $X$ and $Y$, the tensor constraint
$\beta_{X,Y}$ preserves total degree with respect to
these reversed gradings.
\end{lemma}

\begin{proof}

Since $p$ vanishes on $B_r\otimes A_s$ for $r\neq s$,
 {compatibility of the coaction and} \eqref{dual:omega-action}  {give}
\[
 \rho^\Omega(B_r\otimes X_d)\subseteq X_{d-r}.
\]
Let $t_r:X\to B_r\otimes X$ be the  {transpose} used to
construct $\delta_0$ from $p^+$ in
Lemma~\ref{dual:transposition}.
 {Naturality of transposition, applied to the restricted} action
$A_r\otimes X_d\to X_{d+r}$ {,} gives
$t_r(X_d)\subseteq B_r\otimes X_{d+r}$.
 {The finite sums defining} $  \delta_0   $  {therefore restrict to morphisms}
\[
 F_NX\cap X_d  \longrightarrow   
 \bigoplus_{r=0}^{N}B_r\otimes X_{d+r}.
\]
 {These morphisms are compatible with the inclusions for successive}
$  N   $ {, since} $  t_{N+1}   $  {vanishes on} $  F_N   X$ {.
Rationality} and  {exactness of} filtered colimits  {give}
$\varinjlim_N(F_NX\cap X_d)\cong X_d$.
 {Thus taking} the colimit gives a  {factorization}
\[
 X_d\longrightarrow
 \bigoplus_{r\geq0}B_r\otimes X_{d+r}
   \hookrightarrow    B\otimes X
\]
of $\delta_0  |_{X_d}   $.
 {The antipode and the double braiding preserve each}
$B_r\otimes X_{d+r}$ {, by gradedness and naturality}.
 {Equation~}\eqref{dual:omega-coaction}  {consequently gives}
\[
 \delta^\Omega(X_d)\subseteq
 \bigoplus_{r\geq0}B_r\otimes X_{d+r}.
\]
Replacing  {degree} $d$ by $-d$  {makes these} action and coaction
 {formulas compatible with the grading}.
 {For} the  {tensor constraint}
$\beta_{X,Y}=c^{\YD}_{Y,X}c_{Y,X}^{-1}$ {,}
naturality  {of} $  c   $ and \eqref{pre:yd-braiding}  {give}
\[
 \beta_{X,Y}(X_d\otimes Y_e)
 \subseteq
 \bigoplus_{r\geq0}X_{d+r}\otimes Y_{e-r}.
\]
 {Each} summand on the right has   total degree
$-(d+r)-(e-r)=-d-e$ {, equal to that} of $X_d\otimes Y_e$.
 {Thus} $\beta_{X,Y}$ preserves total degree.
\end{proof}

\begin{lemma}\label{refl:lowest}
Let $A=\bigoplus_{n\geq0}A_n$ be a connected graded Hopf algebra in $\D$,
and let $X$ be simple  {in} ${}^{A}_{A}\YD(\D)$.
Suppose that $X=\bigoplus_{d\geq a}X_d$ is a compatible grading
with $X_a\neq0$.
Then $X_a$ is simple in $\D$, the action map $A\otimes X_a\to X$ is an
epimorphism, and
\[
 \prescript{\co A}{}{X}=X_a.
\]
\end{lemma}

\begin{proof}
 {Let} $  \iota_a:X_a\hookrightarrow X   $  {be the inclusion.}
Compatibility of the grading  {implies that} $\delta_X  \iota_a   $
 {factors through} $A_0\otimes X_a$.
 {Connectedness and} the counit identity  {therefore give}
\[
 \delta_X  \iota_a   =\eta_A\otimes  \iota_a   .
\]
Thus $X_a\subseteq\prescript{\co A}{}{X}$.

To prove the reverse direction,
we first show that every nonzero subobject of
$\prescript{\co A}{}{X}$ generates $X$ as an $A$-module.
Let $u:U\hookrightarrow X$ be such a subobject, and  {factor}
\[
 g=\rho_X(\id_A\otimes u),
 \qquad
   A\otimes U\xrightarrow{e}   G(U)  \xrightarrow{k}X   
\]
 {as an epimorphism followed by a monomorphism}.
Associativity of the action makes $G(U)$ an $A$-submodule.
The unit identity gives $U\subseteq G(U)$.
To prove that $G(U)$ is also an $A$-subcomodule,
substitute $\delta_Xu=\eta_A\otimes u$ into the
antipode form of the Yetter--Drinfeld identity
\cite[Remark~3.3]{Simon}.
Naturality of the braiding then gives
\[
 \begin{aligned}
 \delta_Xg
 =(\mu_A\otimes g)
   (\id_A\otimes c_{A\otimes U,A})
   (\id_{A\otimes A\otimes U}\otimes S_A)
   (\Delta_A\otimes c_{A,U})
   (\Delta_A\otimes\id_U).
 \end{aligned}
\]
The right-hand side factors through $  \id_A\otimes k   $ {.
Compose with the quotient map} $A\otimes   X\to A\otimes(X/G(U))   $
 {and cancel the epimorphism} $  e   $ {.
Exactness of} $  A\otimes-   $  {then shows that} $  \delta_Xk   $
 {factors uniquely through} $  \id_A\otimes k   $ {.
This factorization defines the restricted coaction on} $G(U)$.
 {Thus} $G(U)$ is a nonzero Yetter--Drinfeld subobject of $X$ {,
and simplicity} gives $G(U)=X$.

Taking $U=X_a$ proves that the action map
$A\otimes X_a\to X$ is an epimorphism.
 {For any nonzero subobject} $U\subseteq X_a$, the  {same} argument
 {makes} $  A\otimes    U  \to    X$  {a graded epimorphism},  {with} $U$  {placed}
in degree $a$.
 {Its degree-}$  a   $  {component} $A_0\otimes U  \to    X_a$  {is therefore
an epimorphism.
Under} $  A_0\cong\one   $ {, this is the inclusion} $U  \hookrightarrow X_a   $
 {by} the unit identity.  {Hence} $U=X_a$ {, proving simplicity}.

Finally, $\prescript{\co A}{}{X}$ is graded, since
both $\delta_X$ and $\eta_A\otimes\id_X$ preserve degree.
Suppose that it has a nonzero component $U$ of
degree $d>a$.
 {Then} $G(U)=X$,  {whereas} compatibility of the action  {gives}
$
 G(U)\subseteq\bigoplus_{e\geq d}X_e.
$
This contradicts $X_a\neq0$.
 {Hence} $\prescript{\co A}{}{X}=X_a$.

\end{proof}

\subsection{Reflections of tuples}

Let $I=\{1,\ldots,{{r}}\}$, where ${{r}}\geq1$, and let
$M=(M_j)_{j\in I}$ be a tuple of nonzero simple objects of $\C$.
 Let
$
 \B(M):=\B\left(\bigoplus_{j\in I}M_j\right).
$
Fix $i\in I$, and set
\[
 K=\B(M)^{\co\B(M_i)},
\]
where the right coinvariants are taken with respect to the Hopf
projection $\B(M)\to\B(M_i)$ induced by the projection onto $M_i$.
Write $\Omega_i$ for the equivalence of
Theorem~\ref{dual:omega} induced by the pairing
$\B(M_i^*)\otimes\B(M_i)\to\one$ of
Lemma~\ref{dual:pairing}.
For $j\neq i$ and ${{n}}\geq0$, the object $(\ad M_i)^{{n}}(M_j)$ is the
image of the iterated adjoint action
$M_i^{\otimes {{n}}}\otimes M_j\to\B(M)$.
 {Since} $  \B(M_i)   $  {is generated by} $  M_i   $ {, these images give}
\[
 \ad(\B(M_i))(M_j)
 =\bigoplus_{{{n}}\geq0}(\ad M_i)^{{n}}(M_j).
\]

\begin{definition}\label{refl:finite}
The tuple $M$ {admits the $i$-th reflection if}, for every $j\neq i$, there is an
integer $m_{ij}\geq0$ such that
\[
 (\ad M_i)^{m_{ij}}(M_j)\neq0,\qquad
 (\ad M_i)^{m_{ij}+1}(M_j)=0.
\]
Set $L_j:=(\ad M_i)^{m_{ij}}(M_j)$ and define
\[
 R_i(M)_i=M_i^*,\qquad R_i(M)_j=L_j\quad(j\neq i)  .   
\]
We call $R_i(M)$ the $i$-th reflection of $M$.
\end{definition}

 {Each} $  (\ad M_i)^{{n}}(M_j)   $  {belongs to} $  \C   $ {: its defining morphism}
\[
 M_i^{\otimes {{n}}}\otimes M_j\longrightarrow\B(M)({{n}}+1)
\]
 {has source and target} in $\C$ by
Lemma~\ref{pre:nichols-properties}(3) {, and}
$\C  \hookrightarrow\D   $  {is fully faithful and exact}.
If $M$ {admits the $i$-th reflection}, associativity of the adjoint action gives
$(\ad M_i)^{{n}}(M_j)=0$ for ${{n}}>m_{ij}$.
 {The preceding direct sum is then finite, so}
$\ad(\B(M_i))(M_j)$  {and} $  L_j   $  {belong to} $\C$.

\begin{lemma}\label{refl:rational}
\begin{enumerate}
\item
The  {Yetter--Drinfeld module} $K$ is rational if and only if
$M$ {admits the $i$-th reflection}.

\item
Suppose that $M$ {admits the $i$-th reflection}. For each $j\neq i$, the object
\[
 V_j=\Omega_i\bigl(\ad(\B(M_i))(M_j)\bigr)
\]
is simple in ${}_{\B(M_i^*)}^{\B(M_i^*)}\YD(\D)$ {.
Its underlying object} has a compatible positive grading  {defined} by
\[
 (V_j)_{m_{ij}-{{n}}+1}=(\ad M_i)^{{n}}(M_j)
 \qquad(0\leq {{n}}\leq m_{ij}).
\]
 {All other components are zero, and its} degree-one component
$L_j$ is simple in $\C$.
\end{enumerate}
\end{lemma}

\begin{proof}
\textup{(1)}
We write $\cdot$ for the module actions.
On $K$, the action is the restriction of the braided
adjoint action of $\B(M)$.
Suppose first that $M$ {admits the $i$-th reflection}.
 {Associativity and the} grading  {of the adjoint action give}
\[
 \B(M_i)(s)\cdot(\ad M_i)^{{n}}(M_j)
 \subseteq(\ad M_i)^{{{n}}+s}(M_j)
 \qquad({{n}},s\geq0).
\]
 {Thus} $\B(M_i)(s)$ annihilates $\ad(\B(M_i))(M_j)$ for
$s>m_{ij}$, so this module is rational.
Theorem~\ref{coinv:nichols}  {identifies} $K$  {with} the Nichols algebra
of the direct sum of these modules.
 {Closure} under direct sums, tensor products and quotients
 { in Lemma~\ref{dual:rational-closure} therefore makes} $K$ rational.

Conversely, suppose that $K$ is rational.
 {For each} $  j\neq i   $ {, finite presentability of} $M_j$  {and}
$  \varinjlim_NF_NK\cong K   $  {imply that}
$M_j\hookrightarrow K$ factors through some $F_{N_j}K$.
 {Thus} $\B(M_i)({{n}})\cdot M_j=0$  {for} ${{n}}>N_j$.
 { By \eqref{coinv:adjoint-factorization},}
$(\ad M_i)^{{n}}(M_j)  =0   $ for ${{n}}>N_j$.
Since $(\ad M_i)^0(M_j)=M_j\neq0$,  {the nonzero adjoint
components have} a largest  {exponent} $m_{ij}$.
Hence $M$ {admits the $i$-th reflection}.

\textup{(2)}
Fix $j\neq i$.
By Proposition~\ref{coinv:simple},
$\ad(\B(M_i))(M_j)$ is simple in
${}_{\B(M_i)}^{\B(M_i)}\YD(\D)$.
It is rational by \textup{(1)}, so  {the equivalence} $\Omega_i$
sends it to a simple module $V_j$ in
${}_{\B(M_i^*)}^{\B(M_i^*)}\YD(\D)_{\rat}$.
 {Every subobject of} $  V_j   $  {in the  Yetter--Drinfeld category
is} rational by Lemma~\ref{dual:rational-closure}.
 {Hence} $V_j$ is simple in
${}_{\B(M_i^*)}^{\B(M_i^*)}\YD(\D)$.
The multigrading of $K$ gives
\[
 \delta_K\bigl((\ad M_i)^{{n}}(M_j)\bigr)
 \subseteq
 \bigoplus_{s=0}^{{{n}}}
 \B(M_i)(s)\otimes(\ad M_i)^{{{n}}-s}(M_j).
\]
Together with the action formula in \textup{(1)},
this shows that assigning degree ${{n}}$ to
$(\ad M_i)^{{n}}(M_j)$ gives a compatible grading on
$\ad(\B(M_i))(M_j)$.
By Lemma~\ref{dual:grading}, the corresponding component
of $V_j$ has degree $-{{n}}$.
Shifting these degrees by $m_{ij}+1$ gives the asserted
positive grading, with lowest nonzero component
$(V_j)_1=L_j$.

Lemma~\ref{refl:lowest} gives
\[
 \prescript{\co\B(M_i^*)}{}{V_j}=L_j,
 \qquad
 \B(M_i^*)\cdot L_j=V_j,
\]
where the action on $V_j$ is defined by $\Omega_i$.
 {It also} shows that $L_j$ is simple in $\D$.
Since $L_j\in\C$  {and the embedding is exact}, it is simple in $\C$.

\end{proof}

\subsection{Reconstruction from a bosonization}

Let $A$ be a Hopf algebra in $\D$ with bijective antipode, and let
$R$ be a Hopf algebra in
${}^{A}_{A}\YD(\D)$.  { Let} $H=R\#A$  {and write}
\[
 \begin{aligned}
 j&  =\id_R\otimes\eta_A   :R\longrightarrow H,\\
 \iota&  =\eta_R\otimes\id_A   :A\longrightarrow H,\\
 \pi&  =\varepsilon_R\otimes\id_A   :H\longrightarrow A.
 \end{aligned}
\]
\begin{lemma}\label{refl:primitive-splitting}
There is a split exact sequence in $\D$
\[
 0\longrightarrow\prescript{\co A}{}{P(R)}
 \xrightarrow{\ j\ }P(H)
 \xrightarrow{\ \pi\ }P(A)\longrightarrow0.
\]
 {Here the arrows are induced by} $  j   $  {and} $  \pi   $ {, and the section
is induced by} $  \iota   $.
\end{lemma}

\begin{proof}
The Hopf algebra maps $\pi$ and $\iota$ preserve primitives.
Since $\pi\iota=\id_A$, they induce a split epimorphism
$P(H)\to P(A)$ with section induced by $\iota$.
It remains to identify its kernel.

Let $p:P(R)\hookrightarrow R$ be the inclusion.
The bosonization coproduct gives
\begin{equation}\label{refl:primitive-coproduct}
 \Delta_Hjp
 =jp\otimes\eta_H+(\iota\otimes j)\delta_Rp.
\end{equation}
On $\prescript{\co A}{}{P(R)}$, the second term is
$\eta_H\otimes jp$. Thus $j$ maps this subobject into $P(H)$.
Moreover,
$
 \pi j=\eta_A\varepsilon_R,
 \
 \varepsilon_Rp=0.
$
Hence its image lies in the kernel of $P(H)\to P(A)$.

For the reverse inclusion, let
$
 T=\ker\bigl(\pi:P(H)\longrightarrow P(A)\bigr),
$
and let $t:T\hookrightarrow H$ be the composite of the
kernel inclusion $T\hookrightarrow P(H)$ with
$P(H)\hookrightarrow H$.
Then
\[
 \pi t=0,
 \qquad
 \Delta_Ht=t\otimes\eta_H+\eta_H\otimes t.
\]
 {Applying} $  \id_H\otimes\pi   $  {gives}
\[
   (\id_H\otimes\pi)\Delta_Ht=t\otimes\eta_A.   
\]
The right coinvariants of $\pi:H\to A$ are $j(R)$, { 
see \cite[Corollary~3.9.3(1)]{rootsys}}.
Thus there is a unique morphism $v:T\to R$ such that $t=jv$.
 Let $\vartheta=\id_R\otimes\varepsilon_A:H\to R$.
Then $\vartheta j=\id_R$, $\vartheta\eta_H=\eta_R$, and
the bosonization coproduct gives
\[
 (\vartheta\otimes\vartheta)\Delta_Hj=\Delta_R.
\]
Applying $\vartheta\otimes\vartheta$ to the primitive
identity for $t$ yields
\[
 \Delta_Rv
 =(\vartheta\otimes\vartheta)\Delta_Ht
 =v\otimes\eta_R+\eta_R\otimes v.
\]
Thus $v$ factors through $P(R)$.
We may therefore apply \eqref{refl:primitive-coproduct}
to $v$. Comparing it with the primitive identity for
$t=jv$ gives
\[
 (\iota\otimes j)\delta_Rv
 =\eta_H\otimes jv
 =(\iota\otimes j)(\eta_A\otimes v).
\]
Since $\iota\otimes j$ is a monomorphism, we obtain
$\delta_Rv=\eta_A\otimes v$.
Hence $v$ factors through
$\prescript{\co A}{}{P(R)}$.
 {Thus} $  j   $  {identifies} $  \prescript{\co A}{}{P(R)}   $  {with} the kernel.

\end{proof}

Let $N\in\C$, and let
$V=\bigoplus_{j\in J}V_j$ be a finite direct sum of simple objects in
${}^{\B(N)}_{\B(N)}\YD(\D)$.  Suppose each $V_j$ has a compatible grading
supported in positive degrees, with $(V_j)_1\neq0$.  Set
\[
 L=\bigoplus_{j\in J}(V_j)_1,
 \qquad R=\B(V),\qquad H=R\#\B(N).
\]
\begin{lemma}\label{refl:reconstruction}
There is a unique isomorphism of Hopf algebras in $\D$
\[
 \B(L\oplus N)\xrightarrow{\ \sim\ }H
\]
whose restriction to $L\oplus N$ is the canonical inclusion.
\end{lemma}

\begin{proof}
Compatibility of the grading implies that the quantum
symmetrizers preserve total degree. Thus $R$ inherits a grading
$
 R=\bigoplus_{p\geq0}R_p.
$
Here $p$ is the sum of the given degrees of the tensor factors.
The induced action and coaction on $R$ remain compatible
with this grading.
Give $\B(N)$ its standard grading and $H$ the total grading
\[
 H_n=\bigoplus_{p+q=n}R_p\otimes\B(N)(q).
\]
The bosonization formulas make $H$ a graded Hopf algebra.
Since $V$ is supported in positive degrees, we have
$R_0=\one$ and $R_1=L$. Thus
$
 H_0=\one,\ H_1=L\oplus N.
$

Applying Lemma~\ref{refl:lowest} to each $V_j$ gives
\[
 \B(N)\cdot L=V,
 \qquad
 \prescript{\co\B(N)}{}{V}=L.
\]
We now show that $H$ is generated by $H_1$.
The subalgebra of $H$ generated by $L$ and $N$
contains $\B(N)$, since $\B(N)$ is generated by $N$.
The braided adjoint action on $\B(N)\otimes L$ is the composite
\[
\begin{aligned}
 \B(N)\otimes L
 &\xrightarrow{\Delta_{\B(N)}\otimes\id_L}
 \B(N)\otimes\B(N)\otimes L
 \xrightarrow{\id_{\B(N)}\otimes c_{\B(N),L}}
 \B(N)\otimes L\otimes\B(N)\\
 &\xrightarrow{\id_{\B(N)\otimes L}\otimes S_{\B(N)}}
 \B(N)\otimes L\otimes\B(N)
 \xrightarrow{m_H(m_H\otimes\id_H)} H,
\end{aligned}
\]
where the last arrow is the restriction of multiplication
in $H$ along the canonical inclusions.
 {Since this} subalgebra contains $\B(N)$ and $L$, the last
arrow  {factors through it}.  {It therefore contains}
$\ad_H(\B(N))(L)$.
In $H=R\#\B(N)$, the braided adjoint action of $\B(N)$
on $R$ is the composite
\[
 \B(N)\otimes R\xrightarrow{\rho_R}R\hookrightarrow H,
\]
where $\rho_R$ is the given module action, {  see
\cite[Theorem~3.9.2 and Corollary~3.9.3(1)]{rootsys}}.
Restricting to $\B(N)\otimes L$ and using
Lemma~\ref{refl:lowest}, we obtain
\[
 \ad_H(\B(N))(L)=\B(N)\cdot L=V
\]
as subobjects of $H$.
Thus  {this} subalgebra contains $V$,  {and hence} $R$.
 {Multiplication} gives an isomorphism
$
 R\otimes\B(N)\xrightarrow{\ \sim\ }H
$
in $\D$.
 {Since both factors lie in this} subalgebra, it equals $H$.

 { By Lemma~\ref{pre:nichols-properties}(1)},
$P(R)=V$  {and} $P(\B(N))=N$.
 {Since} $\prescript{\co\B(N)}{}{V}=L$,
Lemma~\ref{refl:primitive-splitting}  {identifies} $P(H)$
with $  L\oplus N   $  {through} the  {canonical inclusions}.
Thus $P(H)=H_1$ as subobjects of $H$.
The  {connected graded} Hopf algebra $H$ is generated by $H_1$,
 {so} Lemma~\ref{pre:nichols-characterization} gives the required
isomorphism.
It is unique because $\B(L\oplus N)$ is generated by
$L\oplus N$.

\end{proof}

\begin{theorem}\label{refl:main}
Let $M=(M_j)_{j\in I}$ be a finite tuple of nonzero simple
objects of $\C$, and suppose that $M$ admits the $i$-th reflection.
Then there is a unique isomorphism of Hopf algebras in $\D$
\[
 \Theta_i:\B\bigl(R_i(M)\bigr)
 \xrightarrow{\ \sim\ }
 \Omega_i\!\left(\B(M)^{\co\B(M_i)}\right)\#\B(M_i^*)
\]
whose restriction to
$M_i^*\oplus\bigoplus_{j\neq i}L_j$
is the canonical inclusion.
\end{theorem}

\begin{proof}
By Lemma~\ref{refl:rational},
$K=\B(M)^{\co\B(M_i)}$ is rational.
For each $j\neq i$, the module
\[
 V_j=\Omega_i\bigl(\ad(\B(M_i))(M_j)\bigr)
\]
is simple in ${}_{\B(M_i^*)}^{\B(M_i^*)}\YD(\D)$
and has a compatible positive grading with
$(V_j)_1=L_j$.
Theorem~\ref{coinv:nichols} and
Corollary~\ref{dual:nichols} give an isomorphism
\[
 \B\!\left(\bigoplus_{j\neq i}V_j\right)
   \xrightarrow{\ \sim\ }\Omega_i(K)   
\]
of Hopf algebras in this category  {whose} restriction to  {the
generating object is the inclusion}
$\bigoplus_{j\neq i}V_j  \hookrightarrow\Omega_i(K)   $.
Applying Lemma~\ref{refl:reconstruction} with $N=M_i^*$,
we define $\Theta_i$ as the composite
\[
 {
 \B(R_i(M))
 \xrightarrow{\ \sim\ }
 \B\left(
   \left(\bigoplus_{j\neq i}L_j\right)\oplus M_i^*
 \right)
 \xrightarrow{\ \sim\ }
 \Omega_i(K)\#\B(M_i^*).
 }
\]
The first arrow is induced by the canonical isomorphism
of the generating objects.
This composite has the prescribed restriction to
$M_i^*\oplus\bigoplus_{j\neq i}L_j$.
Uniqueness follows because these subobjects generate
$\B(R_i(M))$.

\end{proof}

\begin{corollary}\label{refl:root-grading}
Under the assumptions of Theorem~\ref{refl:main}, let $N=R_i(M)$.
Then $N$ is a tuple of nonzero simple objects of $\C$ and
admits the $i$-th reflection. Moreover,
\[
 \max\{{{n}}\geq0| (\ad N_i)^{{n}}(N_j)\neq0\}=m_{ij}
 \qquad(j\neq i),
\]
and
\[
 R_i(N)\simeq M.
\]
\end{corollary}

\begin{proof}
The components of $N$ are nonzero and simple by
Lemmas~\ref{refl:rational} and~\ref{pre:simple-dual}.
 Let
\[
 \widetilde{H}=\Omega_i(K)\#\B(M_i^*).
\]
Fix $j\neq i$, and let $V_j$ be as in
Lemma~\ref{refl:rational}.

 {Exactness of} $  \Omega_i   $  {and the canonical} inclusion
$\Omega_i(K)  \hookrightarrow   \widetilde{H}$  {give} a monomorphism
$v_j:V_j\hookrightarrow\widetilde{H}$.
 {For} $  {{n}}\geq0   $ {, let}
$\iota_{j,{{n}}+1}:(V_j)_{{{n}}+1}\hookrightarrow V_j$
 {be} the homogeneous  {inclusion}.
By Lemma~\ref{refl:rational}, the action induces a graded
epimorphism
$
 \B(M_i^*)\otimes L_j\longrightarrow V_j,
$
 {with} $L_j$  {in} degree one. Its component of degree ${{n}}+1$
is an epimorphism
\[
 \rho_{j,{{n}}}:\B(M_i^*)({{n}})\otimes L_j
 \twoheadrightarrow (V_j)_{{{n}}+1}.
\]
Indeed, the  {composite} with the projection onto $(V_j)_{{{n}}+1}$
 {is} an epimorphism  {and} factors through $\rho_{j,{{n}}}$.

Let
$
 \mu_{{n}}:(M_i^*)^{\otimes {{n}}}
 \twoheadrightarrow\B(M_i^*)({{n}})
$
be the multiplication map, with $\mu_0$ the unit isomorphism.
It is an epimorphism because $\B(M_i^*)$ is generated by
$M_i^*$.
Exactness of tensor products implies that
$\mu_{{n}}\otimes\id_{L_j}$ is also an epimorphism.

Let $d_{j,{{n}}}:N_i^{\otimes {{n}}}\otimes N_j  \to   \B(N)$
be the iterated adjoint morphism {, with} $d_{j,0}$ the inclusion
of $N_j$.  {The} Hopf algebra morphism $  \Theta_i   $ preserves
 {adjoint actions and has} the  {prescribed} restrictions to
$N_i=M_i^*$ and $N_j=L_j$ {.
In} $  \widetilde{H}   $ {,} the adjoint action of $\B(M_i^*)$ on
$\Omega_i(K)$ is its module action { , as in the proof of
Lemma~\ref{refl:reconstruction}}.
 {Induction} on ${{n}}$  {and associativity} of  {this} action  {therefore
give} the factorization
\[
 N_i^{\otimes {{n}}}\otimes N_j
 \xrightarrow{\rho_{j,{{n}}}(\mu_{{n}}\otimes\id_{L_j})}
 (V_j)_{{{n}}+1}
 \xrightarrow{\Theta_i^{-1}v_j\iota_{j,{{n}}+1}}
 \B(N)
\]
 {of} $  d_{j,{{n}}}   $ {.
The first arrow} is an epimorphism  {and the second is} a
monomorphism {. Thus} $  \Theta_i^{-1}v_j\iota_{j,{{n}}+1}   $
 {represents the image of} $d_{j,{{n}}}$.
 {Uniqueness} of the image factorization  {and} the grading in
Lemma~\ref{refl:rational}  {give}
\begin{equation}\label{refl:adjoint-reversal}
 (\ad N_i)^{{n}}(N_j)
 \simeq (V_j)_{{{n}}+1}
 \simeq (\ad M_i)^{m_{ij}-{{n}}}(M_j)
 \qquad(0\leq {{n}}\leq m_{ij}).
\end{equation}
For ${{n}}>m_{ij}$,  {one has} $(V_j)_{{{n}}+1}  =0   $, so the same
factorization gives $(\ad N_i)^{{n}}(N_j)=0$.
At ${{n}}=m_{ij}$, \eqref{refl:adjoint-reversal} gives
\[
 (\ad N_i)^{m_{ij}}(N_j)\simeq M_j\neq0.
\]
Thus $N$ admits the $i$-th reflection {,} with the same
 {integers} $  m_{ij}   $.

The adjoint images belong to $\C$  { by the observation after
Definition~\ref{refl:finite}.
Full faithfulness} of $\C\hookrightarrow\D$  {shows that} the
isomorphisms  { in \eqref{refl:adjoint-reversal} lie} in $\C$.
It follows that
\[
 R_i(N)_j\simeq M_j\quad(j\neq i),
 \qquad
 R_i(N)_i=M_i^{**}.
\]
Finally, duality and the braiding give isomorphisms
\[
 \begin{aligned}
 \operatorname{Hom}_{\C}(Y,M_i)
 \cong\operatorname{Hom}_{\C}(M_i^*\otimes Y,\one)
 \cong\operatorname{Hom}_{\C}(Y\otimes M_i^*,\one)
 \cong\operatorname{Hom}_{\C}(Y,M_i^{**}),
 \end{aligned}
\]
natural in $Y\in\C$.
The Yoneda lemma yields $M_i\simeq M_i^{**}$.
Together with the preceding componentwise isomorphisms,
this proves $R_i(N)\simeq M$.

\end{proof}
\begin{remark}\label{refl:semi-cartan}\upshape
Let $M=(M_j)_{j\in I}$ be a finite nonempty tuple of
nonzero simple objects of $\C$.
We say that $M$  {\textbf{admits all reflections}} if
\[
 R_{i_t}\cdots R_{i_1}(M)
\]
is defined for every finite sequence $(i_1,\ldots,i_t)$ in $I$.

Suppose that $M$ admits all reflections.
Let $\mathcal X$ be the set of isomorphism classes of tuples
obtained from $M$ by reflections, including $M$ itself.
For $[P]\in\mathcal X$, define
\[
 \rho_i([P])=[R_i(P)],
 \qquad
 a_{ii}^{[P]}=2,
 \qquad
 a_{ij}^{[P]}
 =-\max\{{{n}}\geq0:(\ad P_i)^{{n}}(P_j)\neq0\}
 \quad(i\neq j).
\]
 Let $A^{[P]}=(a_{ij}^{[P]})_{i,j\in I}$.
These definitions depend only on $[P]$, since isomorphisms
of tuples induce isomorphisms of their Nichols algebras
and adjoint images.

We check that each $A^{[P]}$ is a generalized Cartan matrix.
The diagonal entries are $2$, and the other entries are
nonpositive integers.
For $i\neq j$, let
$
 d_{ij}:P_i\otimes P_j\longrightarrow\B(P)(2)
$
be the adjoint morphism.
The degree-two quantum symmetrizer gives
\[
 \Delta_{1,1}d_{ij}
 =\iota_{ij}
   \bigl(\id_{P_i\otimes P_j}
         -c_{P_j,P_i}c_{P_i,P_j}\bigr),
\]
where $\iota_{ij}$ is the inclusion of $P_i\otimes P_j$
into $(\bigoplus_{k\in I}P_k)^{\otimes2}$.
Since $\Delta_{1,1}$ is a monomorphism on $\B(P)(2)$,
 we obtain
\[
 \begin{aligned}
 a_{ij}^{[P]}=0
 \iff d_{ij}=0
 \iff c_{P_j,P_i}c_{P_i,P_j}
          =\id_{P_i\otimes P_j}
 \iff c_{P_i,P_j}c_{P_j,P_i}
          =\id_{P_j\otimes P_i}
 \iff a_{ji}^{[P]}=0.
 \end{aligned}
\]
Corollary~\ref{refl:root-grading} gives
\[
 \rho_i^2=\id_{\mathcal X},
 \qquad
 a_{ij}^{\rho_i([P])}=a_{ij}^{[P]}
 \quad(i,j\in I).
\]
Thus these data $(I,\mathcal X,(\rho_i)_{i\in I},
  (A^{[P]})_{[P]\in\mathcal X})$ form a semi-Cartan graph in the sense of
\cite[Definition~9.1.1]{rootsys}.
It is connected, since every point of $\mathcal X$
is obtained from $[M]$ by reflections.
\end{remark}

\section{Applications to nonintegral braided fusion categories}
\label{sec:examples}

 {In the categories considered below, the generating object
of a finite Nichols algebra has length at most one.
We prove this for $PSU(2)_k$ at odd levels,
$SU(2)_k$ at levels divisible by four, and the representation
categories of the Virasoro minimal models $\mathcal M(2,p)$ with
$p\ge5$ odd. Meir's results then give consequences for finite
pre-Nichols algebras and finite Hopf algebras.}

Throughout this section, $\C$ is a braided fusion category over
$\mathbb C$. We call a Hopf algebra in $\D$ finite if its underlying
object belongs to $\C$. For $V\in\C$, the algebra $\B(V)$ is finite
if and only if it has only finitely many nonzero homogeneous components.

\subsection{Reflection and multigraded support}
Let $e_1,\dots, e_r$ be the standard basis of $\mathbb{Z}^{r}$ and
let $M=(M_1,\ldots,M_r)$ be a tuple of nonzero simple objects in $\mathcal{C}$ with $\operatorname{deg}M_i=e_i$ for each $i$. Let
\[
 {\operatorname{Supp}\B(M)}=\{\gamma\in\mathbb N_0^r:\B(M)_\gamma\ne0\}.
\]
Thus $\B(M)$ is finite if and only if ${\operatorname{Supp}\B(M)}$ is finite. Suppose
that $\B(M_i)$ is finite, and set
\[
 d_i=\max\{n\mid \B(M_i)(n)\ne0\},\qquad
 m_{ij}=\max\{n\mid (\ad M_i)^n(M_j)\ne0\}\quad(j\ne i).
\]
By \eqref{coinv:adjoint-factorization}, $0\le m_{ij}\le d_i$. Define
$s_i\in\operatorname{GL}(\mathbb Z^r)$ by
\[
 s_i(e_i)=-e_i,\qquad s_i(e_j)=e_j+m_{ij}e_i\quad(j\ne i).
\]
The following consequence of Theorem~\ref{refl:main} will be
used below.
\begin{lemma}\label{app:support}
Suppose that $\B(M_i)$ is finite. Then the map
\begin{equation}\label{app:support-identity}
 {\operatorname{Supp}\B(R_i(M))}\longrightarrow {\operatorname{Supp}\B(M)},\qquad
 \gamma\longmapsto d_i e_i+s_i(\gamma)
\end{equation}
is a bijection. Consequently, $\B(M)$ is finite if and only if
$\B(R_i(M))$ is finite.
\end{lemma}

\begin{proof}
 Let $K_i=\B(M)^{\co\B(M_i)}$, and let ${\operatorname{Supp}K_i}$ be its
multigraded support. Give $\B(M_i^*)(n)$ degree $-ne_i$.
{Retain the original multigrading of $K_i$ on
$\Omega_i(K_i)$. The proof of Lemma~\ref{dual:grading},
applied to each multihomogeneous component, shows that the
transported action sends
$\B(M_i^*)(n)\otimes(K_i)_\gamma$ into
$(K_i)_{\gamma-ne_i}$, while its coaction sends
$(K_i)_\gamma$ into
$\bigoplus_{n\geq0}\B(M_i^*)(n)\otimes(K_i)_{\gamma+ne_i}$.
Moreover, \eqref{dual:tensorator} gives, for
$\gamma,\gamma'\in\mathbb Z^r$,}
\[
 {
 \beta_{K_i,K_i}\bigl((K_i)_\gamma\otimes(K_i)_{\gamma'}\bigr)
 \subseteq\bigoplus_{n\geq0}
 (K_i)_{\gamma+ne_i}\otimes(K_i)_{\gamma'-ne_i}.}
\]
{Thus the action, coaction, and tensor constraint
preserve total multidegree. The transported Hopf algebra
structure and the bosonization therefore preserve it as well.
The canonical generators $M_i^*$ and $L_j$ have degrees
$-e_i=s_i(e_i)$ and $e_j+m_{ij}e_i=s_i(e_j)$, respectively.
Since $\B(R_i(M))$ is generated by these objects,
Theorem~\ref{refl:main} makes $\Theta_i$ graded when its 
component of degree $\gamma$ is assigned degree $s_i(\gamma)$.}

Both $\B(M_i)$ and $\B(M_i^*)$ have nonzero degrees exactly
$0,\ldots,d_i$, by generation in degree one and
Lemma~\ref{dual:pairing}.
Tensor products of nonzero objects in a fusion category
are nonzero. The Radford isomorphism and
Theorem~\ref{refl:main} give
\[
\begin{aligned}
 {\operatorname{Supp}\B(M)}
 &=\{\beta+ne_i\mid
       \beta\in {\operatorname{Supp}K_i},\ n\in\{0,\ldots,d_i\}\},\\
 {\operatorname{Supp}\B(R_i(M))}
 &=\{s_i(\beta-ne_i)\mid
       \beta\in {\operatorname{Supp}K_i},\ n\in\{0,\ldots,d_i\}\}.
\end{aligned}
\]
Since $s_i^2=\id$, the map
$\gamma\mapsto d_i e_i+s_i(\gamma)$ sends the second set
bijectively onto the first, replacing $n$ by $d_i-n$.
This proves both assertions.
\end{proof}

Let $\kappa=(i_1,\ldots,i_l)$ be a reflection sequence of $M$, and write
\[
 R_{(i_1,\ldots,i_k)}(M)=R_{i_k}\cdots R_{i_1}(M),
 \qquad R_{()}(M)=M.
\]
 Let
\[
 \omega= {s_{i_1}\cdots s_{i_l}},
 \qquad
 \beta_k= {s_{i_1}\cdots
 s_{i_{k-1}}({e}_{i_k})}.
\]
Here the $j$-th factor in each product is
$s_{i_j}^{[R_{(i_1,\ldots,i_{j-1})}(M)]}$.
\begin{corollary}\label{app:support-criterion}
Suppose that
$\B\bigl(R_{(i_1,\ldots,i_{k-1})}(M)_{i_k}\bigr)$ is finite,
with highest degree $d_k$, for every $1\le k\le l$. Then the map
\[
 {\operatorname{Supp}\B(R_\kappa(M))}\longrightarrow {\operatorname{Supp}\B(M)},\qquad
 \gamma\longmapsto
 \sum_{k=1}^l d_k\beta_k+\omega(\gamma)
\]
is a bijection. Moreover:
\begin{enumerate}
\item If the matrix of $\omega$ with respect to the standard
basis $(e_i)_{1\leq i \leq r}$ is a negative permutation
matrix, then $\B(M)$ is finite.
\item If $R_\kappa(M)\simeq M$ componentwise and $\omega$
has infinite order as an automorphism of $\mathbb Z^{r}$,
then $\B(M)$ is infinite.
\end{enumerate}
\end{corollary}

\begin{proof}
Iterating Lemma~\ref{app:support} gives the stated bijection.
For \textup{(1)}, write
$\sum_{k=1}^l d_k\beta_k=\sum_{1\leq i \leq r}n_ie_i$.
{Let
$\sum_{i=1}^r a_ie_i\in\operatorname{Supp}\B(M)$.
By surjectivity of the stated bijection, there is
$\gamma\in\operatorname{Supp}\B(R_\kappa(M))$ such that}
\[
 {
 \sum_{i=1}^r a_ie_i
 =\sum_{i=1}^r n_ie_i+\omega(\gamma).
 }
\]
{Since $\gamma\in\mathbb N_0^r$ and $\omega$
is a negative permutation matrix, every coordinate of
$\omega(\gamma)$ is nonpositive. Thus $a_i\le n_i$ for every $i$.
Also, $a_i\ge0$ since
$\operatorname{Supp}\B(M)\subseteq\mathbb N_0^r$. Hence}
\[
 {
 \operatorname{Supp}\B(M)\subseteq
 \left\{
 \sum_{i=1}^r a_ie_i
 \;\middle|\;
 a_i\in\mathbb Z,\quad 0\le a_i\le n_i
 \text{ for every }i
 \right\}.
 }
\]
Hence ${\operatorname{Supp}\B(M)}$ is finite.

For \textup{(2)}, let
\[
 F(\gamma)=\sum_{k=1}^l d_k\beta_k+\omega(\gamma).
\]
Since $R_\kappa(M)\simeq M$ componentwise, their supports
coincide, and $F$ restricts to a bijection from ${\operatorname{Supp}\B(M)}$ to itself.
Suppose that ${\operatorname{Supp}\B(M)}$ is finite. Then there exists $N\ge1$ such
that $F^N(\gamma)=\gamma$ for every $\gamma\in {\operatorname{Supp}\B(M)}$.
The map $F^N$ is affine with linear part $\omega^N$.
Since $0,{e}_1,\ldots,{e}_r\in {\operatorname{Supp}\B(M)}$, we obtain
\[
 \omega^N({e}_i)
 =F^N({e}_i)-F^N(0)
 ={e}_i
 \qquad(1\le i\le r).
\]
Thus $\omega^N=\id$, contrary to the assumption that
$\omega$ has infinite order.
\end{proof}

\begin{proposition}\label{app:rank-two-obstruction}
Let $M=(M_1,M_2)$ and suppose that $\B(M)$ is finite.
Then there is a tuple $P$ obtained from $M$ by finitely many
reflections such that $m_{ij}(P)\le1$ for some distinct
$i,j\in\{1,2\}$. In particular, $(\ad P_i)(P_j)$ is either
zero or simple.
\end{proposition}

\begin{proof}Each $\B(M_i)$ embeds into $\B(M)$ and is therefore finite.
By Lemma~\ref{app:support} and
Corollary~\ref{refl:root-grading}, every tuple obtained from
$M$ by finitely many reflections has simple components and
finite Nichols algebra, and admits both reflections.
Since $\C$ has only finitely many isomorphism classes of
simple objects, there are only finitely many such tuples
up to componentwise isomorphism. As $R_1$ and $R_2$ are
involutive, there exists $n\ge1$ such that
$(R_2R_1)^n(M)\simeq M$ componentwise.
Let $\kappa=(1,2,\ldots,1,2)$ be the corresponding sequence
of length $2n$.

Suppose, to the contrary, that
$m_{12}(P),m_{21}(P)\ge2$ for every tuple $P$ obtained from
$M$ by finitely many reflections.
Let
\[
 \omega= {s_1s_2\cdots s_1s_2}
\]
be the morphism associated with $\kappa$, with
the factors taken along this sequence. With respect to
$(e_1,e_2)$, each consecutive pair of factors has the form
\[
 T(a,b)=
 \begin{pmatrix}-1&a\\0&1\end{pmatrix}
 \begin{pmatrix}1&0\\b&-1\end{pmatrix}
 =
 \begin{pmatrix}ab-1&-a\\b&-1\end{pmatrix},
 \qquad a,b\ge2.
\]
For $x>y\ge0$, write
$T(a,b)(xe_1+ye_2)=x'e_1+y'e_2$. Then
\[
 y'=bx-y\ge2x-y>x,\qquad
 x'=ay'-x\ge2y'-x>y'.
\]
Thus each such matrix preserves the region $x>y\ge0$ and
strictly increases the second coordinate. The same holds for
their product $\omega$. Since $e_1$ lies in this region, the
vectors $\omega^k(e_1)$, $k\ge0$, have strictly increasing
second coordinates. Hence $\omega$ has infinite order.
Together with $R_\kappa(M)\simeq M$, this contradicts
Corollary~\ref{app:support-criterion}\textup{(2)}.

Therefore $m_{ij}(P)\le1$ for some reflected tuple $P$ and
distinct indices $i,j$. If $m_{ij}(P)=0$, then
$(\ad P_i)(P_j)=0$. If $m_{ij}(P)=1$, this first adjoint is
the terminal adjoint and is simple by
Lemma~\ref{refl:rational}.
\end{proof}

\begingroup
\subsection{Nichols algebras in \texorpdfstring{$PSU(2)_k$}{PSU(2)k} at odd levels}
The categories $PSU(2)_3$, $PSU(2)_5$, and $PSU(2)_7$
occur in the classification of unitary modular categories
with at most four simple objects
\cite[Sections~5.3--5.4]{RSW09}.
In particular, $PSU(2)_3$ is the Fibonacci category,
whose braid group representations give a model for universal
topological quantum comletation
\cite[Theorem~6.7]{RSW09}.
The theorem below gives a uniform restriction on finite
Nichols algebras throughout the odd-level family.

Let $k\ge3$ be odd. The modular fusion category $PSU(2)_k=(A_1,k)_{1/2}$
is the full tensor subcategory of $SU(2)_k$ generated by the simple
objects with even labels.  Let $m=(k-1)/2$, and write $V_a$ for the object
with label $2a$, where $0\le a\le m$ and $V_0=\one$.
The fusion rules are
\begin{equation}\label{app:psu-fusion}
 V_a\otimes V_b\simeq
 \bigoplus_{c=|a-b|}^{\min(a+b,k-a-b)}V_c
 \qquad(0\le a,b\le m);
\end{equation}
see \cite[(1.1)]{RS11}.
The twist $\theta$ acts on $V_a$ by the scalar
\begin{equation}\label{app:psu-twists}
 \theta_a=\zeta^{a(a+1)},
\end{equation}
where $\zeta=e^{2\pi i/(k+2)}$ or its inverse, according to the
braiding convention, see  {\cite[Example~5.1, (5.5)]{GMR24}}.

\begin{theorem}\label{app:psu}
Let $k\ge3$ be odd and $V\in PSU(2)_k$.
If $\B(V)$ is finite, then $V$ has length at most one.
\end{theorem}

\begin{proof}
For $1\le a,b\le m$, the balancing identity shows that
$c_{V_b,V_a}c_{V_a,V_b}$ acts on the summand $V_c$ in
\eqref{app:psu-fusion} by $\theta_c/(\theta_a\theta_b)$.
Remark~\ref{refl:semi-cartan} gives
\begin{equation}\label{app:psu-first-adjoint}
 (\ad V_a)(V_b)\simeq
 \bigoplus_{\substack{|a-b|\le c\le\min(a+b,k-a-b)\\
                     \theta_c\ne\theta_a\theta_b}}V_c.
\end{equation}
The object is formed in $\B(V_a\oplus V_b)$, with two
separate summands when $a=b$.

By \eqref{app:psu-first-adjoint}, it suffices to find at least
two summands $V_c$ of $V_a\otimes V_b$ satisfying
$\theta_c\ne\theta_a\theta_b$.
The fusion rule and this condition are symmetric in $a,b$,
so assume $a\le b$.
The number of summands in \eqref{app:psu-fusion} is
\[
 \min(a+b,k-a-b)-(b-a)+1
 =\min(2a+1,k-2b+1)\ge2.
\]
Suppose first that there are at least three  summands,
and choose three consecutive labels.
For any two of these labels $c<d$, we have $d-c\in\{1,2\}$.
If $\theta_c=\theta_d$, then
\[
 (k+2)\mid d(d+1)-c(c+1)
          =(d-c)(c+d+1).
\]
Since $k+2$ is odd, it is coprime to $d-c$.
Hence $(k+2)\mid(c+d+1)$, contrary to
$0<c+d+1<k+2$.
Thus the three twists are pairwise distinct.
At most one equals $\theta_a\theta_b$, so at least two
of the chosen summands occur in
\eqref{app:psu-first-adjoint}.

Suppose now that there are exactly two fusion summands.
Since $2a+1\ge3$, the count above gives $k-2b+1=2$,
and hence $b=m$.
The fusion rule becomes
\[
 V_a\otimes V_m\simeq V_{m-a}\oplus V_{m-a+1}.
\]
Using $k+2=2m+3$ in \eqref{app:psu-twists}, we obtain
\[
 \frac{\theta_{m-a}}{\theta_a\theta_m}=\zeta^a,
 \qquad
 \frac{\theta_{m-a+1}}{\theta_a\theta_m}=\zeta^{-a-1}.
\]
Neither ratio is $1$, since $\zeta$ has order $2m+3$
and $1\le a\le m$.
Thus both summands occur in \eqref{app:psu-first-adjoint}.

In all cases, $(\ad V_a)(V_b)$ has at least two simple
summands and is therefore nonzero and not simple.
Suppose now that $\B(V)$ is finite and $V$ has length at least two.
Choose simple objects $M_1,M_2$ such that $M_1\oplus M_2$
is a direct summand of $V$.
We have an inclusion
$\B(M_1\oplus M_2)\hookrightarrow\B(V)$, so
$\B(M_1\oplus M_2)$ is finite.
Proposition~\ref{app:rank-two-obstruction} gives a reflected tuple
$P=(P_1,P_2)$ for which $(\operatorname{ad}(P_1)(P_2))$ is zero or simple.
Its Nichols algebra is finite by Lemma~\ref{app:support}.
Since $\B(P_i)$ embeds into $\B(P)$ and $\B(\one)$ is infinite,
neither $P_i$ is isomorphic to $\one$.
This contradicts the preceding calculation.
\end{proof}

\begingroup
\subsection{Nichols algebras in \texorpdfstring{$SU(2)_{4m}$}{SU(2)4m}}
The categories $SU(2)_k$ arise by semisimplifying tilting
modules for $U_q(\mathfrak{sl}_2)$ at roots of unity
\cite[Section~5.1]{GMR24}.
Their modular structures give rise to three-dimensional
topological quantum field theories
\cite[Section~1]{RSW09}.
We study finite Nichols algebras in these categories at
levels divisible by four.

Let $k\ge4$ be divisible by $4$. Write $V_0,\ldots,V_k$ for
the simple objects of $SU(2)_k$, with $V_0=\one$.  Let
$n=k+2$ and $\zeta=e^{2\pi i/n}$. The fusion rules and twists are
\begin{equation}\label{app:su2-data}
 V_a\otimes V_b\simeq
 \bigoplus_{\substack{|a-b|\le c\le\min(a+b,2k-a-b)\\
                      c\equiv a+b\pmod2}}V_c,
 \qquad
 \theta_a=\exp\!\left(\frac{\pi i a(a+2)}{2n}\right);
\end{equation}
see \cite[Example~5.1]{GMR24}.
In particular, $V_1$ has Frobenius--Perron dimension
$2\cos(\pi/(k+2))\notin\mathbb Z$. So $SU(2)_k$ is not an integral fusion category.

\begin{theorem}\label{app:su2}
Let $k\ge4$ be divisible by $4$ and $V\in SU(2)_k$.
If $\B(V)$ is finite, then $V$ has length at most one.
\end{theorem}

\begin{proof}
The invertible simple objects are $V_0$ and $V_k$.
 {On an invertible simple object, the self-braiding is
the scalar given by the quadratic form of the pointed subcategory.
For $V_k$, this scalar is $\theta_k=e^{\pi i k/2}=1$, since $4\mid k$,
see \cite[Sections~3.1--3.2 and Example~5.1, (5.5)]{GMR24}. Hence}
\[
  {c_{V_k,V_k}=\id_{V_k\otimes V_k}.}
\]
 {Its degree-$d$ quantum symmetrizer is therefore
$d!\,\id$, so $\B(V_k)$ is infinite. The same holds for
$\B(V_0)=\B(\one)$.}
As in the last paragraph of the proof of Theorem~\ref{app:psu},
Proposition~\ref{app:rank-two-obstruction} reduces the assertion
to showing that $(\ad V_a)(V_b)$ has length at least two whenever
$1\le a,b\le k-1$.

By  {the balancing identity} and Remark~\ref{refl:semi-cartan}, the object  {$(\ad V_a)(V_b)$} is isomorphic to the direct sum of those summands $V_c$ in
\eqref{app:su2-data} for which $\theta_c\ne\theta_a\theta_b$.
As before, the adjoint is formed using two separate summands
when $a=b$.
Assume $a\le b$. There are $\min(a,k-b)+1$  summands,
so at least two.

If there are at least four summands, take four consecutive
labels. Adjacent labels $c,c+2$ have distinct twists, since
\[
 \frac{\theta_{c+2}}{\theta_c}=\zeta^{c+2}\ne1
 \qquad(0\le c\le k-2).
\]
At least one summand from each of the first and last pairs
has twist different from $\theta_a\theta_b$.
Thus at least two summands occur in the first adjoint.

For the remaining cases, note that
$\theta_{k-c}=(-1)^c\theta_c$.
If $V_c$ occurs in $V_a\otimes V_d$, then $V_{k-c}$ occurs
in $V_a\otimes V_{k-d}$ and $c\equiv a+d\pmod2$. Hence
\begin{equation}\label{app:su2-twist-ratio}
 \frac{\theta_{k-c}}{\theta_a\theta_{k-d}}
 =(-1)^a\frac{\theta_c}{\theta_a\theta_d}.
\end{equation}
If there are exactly two summands, then $a=1$ or $b=k-1$.
For $V_1\otimes V_b$, the two ratios
$\theta_c/(\theta_1\theta_b)$ are
\[
 \exp\!\left(\frac{\pi i b}{n}\right),
 \qquad
 \exp\!\left(-\frac{\pi i(b+2)}{n}\right).
\]
Neither is $1$ or $-1$, since $1\le b\le k-1$.
Formula~\eqref{app:su2-twist-ratio} gives the same conclusion
for $V_a\otimes V_{k-1}$. Thus both summands occur in the
first adjoint.

If there are exactly three summands, then $a=2$ or $b=k-2$.
For $V_2\otimes V_b$, where $2\le b\le k-2$, the ratios
$\theta_c/(\theta_2\theta_b)$ are
\[
 \zeta^b,\qquad \zeta^{-2},\qquad \zeta^{-b-2}.
\]
None is $1$. On $V_a\otimes V_{k-2}$, formula
\eqref{app:su2-twist-ratio} gives
\[
 (-1)^a\zeta^a,\qquad
 (-1)^a\zeta^{-2},\qquad
 (-1)^a\zeta^{-a-2},
 \qquad 2\le a\le k-2.
\]
For even $a$, none is $1$. For odd $a$, the middle scalar
is not $1$, since $n\ge6$. If both outer scalars were $1$,
their product would give $\zeta^{-2}=1$, again impossible.
Thus at least two summands occur in every case.

If a finite Nichols algebra had at least two simple generating
summands, Proposition~\ref{app:rank-two-obstruction}, applied
to two of them, would give a reflected pair with a zero or
simple first adjoint.
Both components would have finite Nichols algebras and hence
be noninvertible, contradicting the preceding calculation.
\end{proof}

\subsection{Virasoro minimal models}
\label{app:minimal-models}
{
Virasoro minimal models are basic examples of rational
conformal field theory, and their representation categories
are modular tensor categories
\cite[Section~2.2]{MS24}.
The family $\mathcal M(2,p)$ considered below includes the
Yang--Lee model and provides nonunitary examples
\cite[Section~5.2]{BD12}.
We show that finite Nichols algebras in these representation
categories satisfy the same restriction on their generating
objects.
}

Let $p\ge5$ be odd. The Virasoro minimal model $\mathcal M(2,p)$
is the simple Virasoro vertex operator algebra of central charge
$1-3(p-2)^2/p$. Its finite direct sums of irreducible modules form
a modular tensor category, denoted by
$\operatorname{Rep}\mathcal M(2,p)$, see
\cite[Section~2.2]{MS24}.
The irreducible modules have labels $L_{1,s}$, $1\le s\le p-1$,
with $L_{1,s}\simeq L_{1,p-s}$. We choose the representatives
\[
 X_a=L_{1,2a+1},\qquad 0\le a\le (p-3)/2.
\]
Their fusion rules and conformal weights give
\begin{align}
 X_a\otimes X_b
 &\simeq\bigoplus_{c=|a-b|}^{\min(a+b,p-2-a-b)}X_c,
 \label{app:minimal-fusion}\\
 h_{1,2a+1}&=\frac{2a(a+1)}p-a,
 \qquad
 \theta_{X_a}=e^{2\pi i h_{1,2a+1}}
             =\zeta^{2a(a+1)},\qquad \zeta=e^{2\pi i/p}.
 \label{app:minimal-twist}
\end{align}
Here  $X_0=\one$,
see \cite[Section~5.1]{BD12}.

\begin{corollary}\label{app:minimal}
Let $p\ge5$ be odd and
$V\in\operatorname{Rep}\mathcal M(2,p)$.
If $\B(V)$ is finite, then $V$ has length at most one.
\end{corollary}

\begin{proof}
The fusion rules \eqref{app:minimal-fusion} coincide with
\eqref{app:psu-fusion} for $k=p-2$.
The twists \eqref{app:minimal-twist} are obtained from
\eqref{app:psu-twists} by replacing $\zeta$ with $\zeta^2$.
Since $p$ is odd, $\zeta^2$ has order $p$.
In particular, for a summand $X_c$ of $X_a\otimes X_b$,
\[
 \theta_{X_c}=\theta_{X_a}\theta_{X_b}
 \quad\Longleftrightarrow\quad
 p\mid c(c+1)-a(a+1)-b(b+1).
\]
This is the same condition as in the proof of
Theorem~\ref{app:psu}.
The counting argument there therefore shows that
$(\ad X_a)(X_b)$ has at least two simple summands
whenever $X_a$ and $X_b$ are nonunit simple objects.

{
If a finite Nichols algebra had at least two simple generating
summands, Proposition~\ref{app:rank-two-obstruction}, applied
to two of them, would give a reflected pair with a zero or
simple first adjoint.
Its components cannot be the unit object, since their Nichols
algebras are finite whereas $\B(\one)$ is infinite.
This contradicts the preceding calculation.
}
\end{proof}

These categories are nonintegral. Indeed, their fusion rules give
\[
 \operatorname{FPdim}(X_a)
 =\frac{\sin((2a+1)\pi/p)}{\sin(\pi/p)},
 \qquad
 \operatorname{FPdim}(X_{(p-3)/2})=2\cos(\pi/p)\in(1,2).
\]
For $p=5$, this is the Yang--Lee category, with
$X_1\otimes X_1\simeq\one\oplus X_1$ and
$\theta_{X_1}=e^{-2\pi i/5}$.
It is the nonunitary Galois conjugate of the Fibonacci category
$PSU(2)_3$, see \cite[Section~1]{RSW09} and
\cite[Section~5.2]{BD12}.
Corollary~\ref{app:minimal} implies that a finite Nichols algebra
in this category is either $\one$ or generated by $X_1$.

\endgroup

\subsection{Finite Hopf algebras}
The generation problem for finite-dimensional pointed Hopf
 {algebras leads to a  question for the associated
graded Hopf algebra of the coradical filtration: is its braided
Hopf algebra of coinvariants generated by its primitive subobject?}
Meir relates this question to rigidity of finite Nichols
algebras under deformation
\cite[Section~1 and Theorems~1.6--1.7]{defofnichols}.
The length restrictions proved above allow us to apply his
results to finite pre-Nichols algebras and to finite Hopf
algebras whose semisimple quotient and coradical are both
the unit object.

A pre-Nichols algebra is a connected graded Hopf algebra generated
by its degree-one component. For a Hopf algebra $H\in\C$, write
$J(H)$ for its Jacobson radical and $\operatorname{corad}(H)$ for
its coradical. These are the largest nilpotent ideal and the
largest cosemisimple subcoalgebra, respectively, see
\cite[Section~2.2]{defofnichols}.

\begin{corollary}\label{app:psu-hopf}
{Let $\C$ be one of the categories $PSU(2)_k$ for odd
$k\ge3$, $SU(2)_k$ for $4\mid k$ and $k\ge4$, or
$\operatorname{Rep}\mathcal M(2,p)$ for odd $p\ge5$.}
\begin{enumerate}
\item Every finite pre-Nichols algebra in $\C$ is a Nichols algebra.
\item Let $H$ be a Hopf algebra in $\C$ such that
\[
 H/J(H)\simeq\one,
 \qquad \operatorname{corad}(H)\simeq\one.
\]
Then $H\simeq\B(P(H))$ as Hopf algebras.
Moreover, $P(H)$ is zero or simple, and $P(H)\not\simeq\one$.
\end{enumerate}
\end{corollary}

\begin{proof}
By  {Theorems~\ref{app:psu} and~\ref{app:su2}, and
Corollary~\ref{app:minimal}}, every finite Nichols algebra in $\C$
is either $\one$ or generated by a simple object.
By \cite[Theorem~1.4]{defofnichols}, those generated by simple
objects are rigid in the sense of \cite[Definition~1.3]{defofnichols}.
The unit Hopf algebra is also rigid, since every Hopf algebra
structure on its underlying object is isomorphic to it.
Thus \cite[Theorem~1.6]{defofnichols} gives \textup{(1)} and the
isomorphism in \textup{(2)}.

Applying  {the same length restriction} to $\B(P(H))\simeq H$ shows
that $P(H)$ is zero or simple. It cannot be isomorphic to
$\one$, since $\B(\one)$ is infinite.
\end{proof}

\endgroup

\section*{ {Declaration of AI assistance}}
 {The author used ChatGPT to improve the exposition,
and check arguments. The author  takes full responsibility for the manuscript.}

\bibliographystyle{plain}
      \bibliography{ref}
    \end{document}